\documentclass{gtart}
\usepackage{amssymb}
\usepackage{amsmath,amscd}

\input gtmonout
\volumenumber{2}
\volumeyear{1999}
\volumename{Proceedings of the Kirbyfest}
\pagenumbers{349}{406}
\papernumber{20}
\received{30 December 1998}\revised{29 March 1999}
\published{21 November 1999}

\def\S{section \ignorespaces}

\newtheorem{theorem}{Theorem}[section]
\newtheorem{lemma}[theorem]{Lemma}
\newtheorem{proposition}[theorem]{Proposition}

\theoremstyle{definition}
\newtheorem{definition}[theorem]{Definition}
\newtheorem{problem}[theorem]{Problem}
\newtheorem{conjecture}[theorem]{Conjecture}
\newtheorem*{remark}{Remark}

\def\varph{{\varphi}}
\def\ra{{\rightarrow}}
\def\lra{{\longrightarrow}}
\def\om{{\omega}}

\def\la{{\lambda}}
\def\al{{\alpha}}
\def\La{{\Lambda}}
\def\ga{{\gamma}}
\def\Ga{{\Gamma}}
\def\bC{{\mathbb C}}
\def\bZ{{\mathbb Z}}
\def\bQ{{\mathbb Q}}
\def\bR{{\mathbb R}}
\def\inv{{^{-1}}}
\def\Sg{{\Sigma}_g}
\def\Mg{{\cal M}_g}
\def\Ng{{\cal N}_g}

\def\Tg{{\cal T}_g}
\def\hg{{\mathfrak h}_g}
\def\tg{{\mathfrak t}_g}

\def\ng{{\mathfrak n}_g}
\def\Ig{{\cal I}_g}
\def\Kg{{\cal K}_g}
\def\Msg{{\cal M}_{g,*}}
\def\Nsg{{\cal N}_{g,*}}

\def\Tsg{{\cal T}_{g,*}}
\def\tsg{{\mathfrak t}_{g,*}}

\def\nsg{{\mathfrak n}_{g,*}}
\def\Isg{{\cal I}_{g,*}}

\def\hsg{{\mathfrak h}_{g,*}}
\def\M1g{{\cal M}_{g,1}}
\def\N1g{{\cal N}_{g,1}}
\def\U1g{{\cal U}_{g,1}}
\def\T1g{{\cal T}_{g,1}}
\def\h1g{{\mathfrak h}_{g,1}}
\def\n1g{{\mathfrak n}_{g,1}}
\def\I1g{{\cal I}_{g,1}}
\def\K1g{{\cal K}_{g,1}}
\def\vGa{{\varGamma}}

\newcommand\Hom{\operatorname{Hom}}

\newcommand\Aut{\operatorname{Aut}}

\newcommand\Out{\operatorname{Out}}

\newcommand\sgn{\operatorname{sgn}}

\newcommand\Gr{\operatorname{Gr}}

\newcommand\ho{\operatorname{Homeo}}
\newcommand\Gal{\operatorname{Gal}}
\newcommand\Ker{\operatorname{Ker}}
\newcommand\Cok{\operatorname{Cok}}
\newcommand\im{\operatorname{Im}}
\newcommand\Int{\operatorname{Int}}

\newcommand\Tr{\operatorname{Tr}}

\newcommand\BD{\operatorname{BDiff}}
\newcommand\ED{\operatorname{EDiff}}

\newcommand\Diff{\operatorname{Diff}}

\begin{document}

\shorttitle{Structure of the mapping class groups of surfaces}
\title{Structure of the mapping class groups of surfaces: \\
a survey and a prospect}

\author{Shigeyuki Morita}
\address{Department of Mathematical Sciences, 
University of Tokyo \\Komaba, Tokyo 153-8914, Japan}
\email{morita@ms.u-tokyo.ac.jp}

\begin{abstract}
In this paper, we survey recent works on the structure of the
mapping class groups of surfaces
mainly from the point of view of topology.
We then discuss several possible directions
for future research. These include the relation between
the structure of the mapping class group and invariants of 3--manifolds,
the unstable cohomology of the moduli space of curves and
Faber's conjecture,
cokernel of the Johnson homomorphisms and the Galois 
as well as other new obstructions,
cohomology of certain infinite dimensional
Lie algebra and
characteristic classes of outer automorphism groups
of free groups and the secondary characteristic 
classes of surface bundles.
We give some experimental results concerning each of them
and, partly based on them, 
we formulate several conjectures and problems.
\end{abstract}


\maketitlepage

\cl{\small\it This paper is dedicated to Robion C Kirby 
on the occasion of his $60^{th}$ birthday.}

\section{Introduction}

Let $\Sg$ be a closed oriented surface of genus $g\geq 2$
and let $\Mg$ be its
mapping class group. This is the group 
consisting of path components of $\Diff_+\Sg$, 
which is the group of orientation
preserving diffeomorphisms of $\Sg$. 
$\Mg$ acts on the Teichm\"uller space
${\cal T}_g$ of $\Sg$ properly discontinuously and the quotient space
${\mathbf M}_g = {\cal T}_g/\Mg$
is the (coarse) moduli space of curves of genus $g$.
${\cal T}_g$ is known to be homeomorphic
to $\bR^{6g-6}$.
Hence we have a natural isomorphism
$$
H^*(\Mg;\bQ) \cong H^*({\mathbf M}_g;\bQ).
$$
On the other hand, 
by a theorem of Earle--Eells \cite{EE}, the identity component of 
$\Diff_+\Sg$ is contractible for $g\geq 2$ 
so that the classifying space $\BD_+\Sg$ is an Eilenberg--MacLane
space $K(\Mg,1)$. Therefore we have also a natural isomorphism
$$
H^*(\BD_+\Sg)\cong H^*(\Mg).
$$
Thus the mapping class group serves as the orbifold fundamental 
group of the moduli space ${\mathbf M}_g$ and at the same time
it plays the role of the {\it universal monodromy group} for
oriented $\Sg$--bundles.
Any cohomology class of the mapping class group
can be considered as a characteristic class of oriented
surface bundles and, over the rationals, it can also be
identified as a cohomology class of the moduli space.
\par
The Teichm\"uller space ${\cal T}_g$ and the moduli space ${\mathbf M}_g$
are important
objects primarily in complex analysis and 
algebraic geometry.
Many important results 
concerning these two spaces have been obtained
following the fundamental works of Ahlfors, Bers  
and Mumford.
Because of the limitation of our knowledge, we only
mention here a survey paper of Hain and Looijenga \cite{HL}
for recent works on ${\mathbf M}_g$,
mainly from the viewpoint of algebraic geometry, 
and a book by Harris and Morrison \cite{HM} 
for basic facts as well as more advanced results.
\par
From a topological point of view, 
fundamental works of Harer \cite{Har1, Har2} on the homology of the
mapping class group and also of Johnson (see \cite{Jo4}) on the structure
of the Torelli group,
both in early $80$'s, 
paved the way towards modern topological studies of 
$\Mg$ and ${\mathbf M}_g$.
Here the Torelli group, denoted by $\Ig$, is
the subgroup of $\Mg$ consisting of those
elements which act on the homology of $\Sg$ trivially.
\par
Slightly later, the author began a study 
of the classifying space 
$\BD_+\Sg$ of surface bundles which also belongs to topology.
The intimate relationship between three universal spaces,
${\cal T}_g, {\mathbf M}_g$ and $\BD_+\Sg$ described above,
imply that there should exist various interactions 
among the studies of these spaces which are peculiar to 
various branches of mathematics including the ones mentioned above. 
Although it is not always easy to understand mutual viewpoints,
we believe that doing so will enhance individual understanding
of these spaces.
\par
In this paper, we would like to survey some aspects of recent
topological study of the mapping
class group as well as the moduli space.
More precisely, we focus on a study of the mapping class
group which is related to the structure of the Torelli group
$\Ig$ together with a natural action of the Siegel modular
group $Sp(2g,\bZ)$ on some graded modules associated with
the lower (as well as other) central series of $\Ig$.
Here it turns out that 
explicit descriptions of $Sp$--invariant tensors of various
$Sp$--modules using
classical symplectic representation theory,
along the lines of Kontsevich's fundamental works in \cite{Ko2, Ko3},
and also
Hain's recent work \cite{H3} on $\Ig$ using mixed Hodge structures
can play very important roles.
These two points will be reviewed in \S4 and \S 5, respectively.
In the final section (\S 6), we describe several experimental
results, with sketches of proofs, by which we would like to
propose some possible directions for future research.
\par
This article can be considered as a continuation of our earlier
papers \cite{Mrt7, Mrt8, Mrt11}.

\rk{Acknowledgements}
We would like to express our hearty thanks to R Hain, N Kawazumi and
H Nakamura for many enlightening discussions and helpful information.
We also would like to thank C Faber, S Garoufalidis, M Kontsevich,
J Levine, E Looijenga,
M Matsumoto, J Murakami and K Vogtmann for helpful discussions 
and communications.
Some of the explicit computations described in \S 6 were done 
by using Mathematica. It is a pleasure to thank M Shishikura for help
in handling Mathematica.

\section{$\Mg$ as an extension of the Siegel modular group by the Torelli group}

Let us simply write $H$ for $H_1(\Sg,\bZ)$. We have the intersection
pairing
\begin{equation*}
\mu\co H\otimes H\lra \bZ
\end{equation*}
which is a non-degenerate skew symmetric bilinear form on $H$. 
The natural action
of $\Mg$ on $H$, which preserves this pairing, induces the
classical representation
\begin{equation*}
\rho_0\co \Mg\lra \Aut H.
\end{equation*}
If we fix a symplectic basis of $H$, then $\Aut H$ can be
identified with the Siegel modular group $Sp(2g,\bZ)$ so
that we can write
\begin{equation*}
\rho_0\co \Mg\lra Sp(2g,\bZ).
\end{equation*}
The {\it Torelli group}, denoted by $\Ig$, is defined to be the
kernel of $\rho_0$. Thus we have the following basic extension
of three important groups
\begin{equation}
1\lra \Ig\lra\Mg\lra Sp(2g,\bZ)\lra 1.
\label{1}
\end{equation}

Associated to each of these groups, we have various moduli spaces.
Namely the (coarse) moduli space 
${\mathbf M}_g$ of genus $g$ curves for $\Mg$,
the moduli space ${\mathbf A}_g$ of principally polarized abelian
varieties for $Sp(2g,\bZ)$ and the {\it Torelli space} ${\mathbf T}_g$
for $\Ig$. Here the Torelli space is defined to be the
quotient of the Teichm\"uller space ${\cal T}_g$ by the natural
action of $\Ig$ on it.
Since $\Ig$ is known to be torsion free, 
${\mathbf T}_g$ is a complex manifold.
We have holomorphic mappings between these moduli spaces
$$
{\mathbf T}_g\lra {\mathbf M}_g\lra {\mathbf A}_g
$$
where the first map is an infinite ramified covering and the
second map is injective by the theorem of Torelli.
\par
By virtue of the above facts, we can investigate the structure of $\Mg$
(or that of ${\mathbf M}_g$)
by combining individual study of $\Ig$ and $Sp(2g,\bZ)$ 
(or ${\mathbf T}_g$ and ${\mathbf A}_g$)
together with some additional investigation of the action 
of $Sp(2g,\bZ)$ on the structure of the Torelli group or Torelli space.
Here it turns out that the symplectic representation theory
can play a crucial role.
However, before reviewing them, let us first recall the 
fundamental works of
D Johnson on the structure of $\Ig$ very briefly (see \cite{Jo4} for
details)
because it is the starting point of the above method.
\par
Johnson proved in \cite{Jo3} that $\Ig$ is finitely 
generated for all $g\geq 3$ by constructing explicit generators
for it. Before this work, a homomorphism
$$
\tau\co \Ig\lra \La^3 H/H
$$
was introduced in \cite{Jo2} which generalized an earlier work of
Sullivan \cite{Su1} extensively and is now called the Johnson homomorphism.
Here $\La^3 H$ denotes the third exterior power of $H$ and
$H$ is considered as a natural submodule of $\La^3 H$
by the injection
$$
H\ni u\longmapsto u\land \om_0\in \La^3 H
$$
where $\om_0\in \La^2 H$ is the symplectic class (in homology)
defined as $\om_0=\sum_i x_i\land y_i$ for any symplectic basis
$x_i, y_i\ (i=1,\cdots, g)$ of $H$. 
\par
Let
${\cal K}_g\subset \Mg$ be the subgroup of $\Mg$ generated by
all Dehn twists along {\it separating} simple closed curves on
$\Sg$. It is a normal subgroup of $\Mg$ and is contained in the Torelli
group $\Ig$. In \cite{Jo5}, Johnson proved that $\Kg$ is exactly equal 
to $\Ker\tau$ so that we have an exact sequence
\begin{equation}
1\lra \Kg\lra \Ig\overset{\tau}{\lra} \La^3 H/H\lra 1.
\label{2}
\end{equation}

Finally in \cite{Jo6}, he determined the abelianization
of $\Ig$ for $g\geq 3$ in terms of certain combination of $\tau$ 
and the totality of the Birman--Craggs homomorphisms defined in \cite{BC}. 
The target of the latter homomorphisms are $\bZ/2$
so that the
first rational homology group of $\Ig$
(or more precisely, the abelianization of $\Ig$ modulo 2 torsions)  
is given simply by $\tau$. Namely we have an isomorphism
$$
\tau\co H_1(\Ig;\bQ)\cong \La^3 H_\bQ/H_\bQ
$$
where $H_\bQ=H\otimes\bQ$.
\par

\begin{problem}
Determine whether the Torelli group $\Ig\ (g\geq 3)$ is
finitely presentable or not.
If the answer is yes, give an explicit finite presentation of it.
\label{prob:torelli}
\end{problem}

It should be mentioned here that Hain \cite{H3} 
proved that the Torelli Lie algebra ${\mathfrak t}_g$, 
which is the Malcev Lie algebra
of $\Ig$, is finitely presentable for all $g\geq 3$.
Moreover he gave an explicit finite presentation of ${\mathfrak t}_g$
for any $g \geq 6$ which turns out to be very simple,
namely there arise only quadratic relations.
Here a result of Kabanov \cite{Ka} played an important role 
in bounding the degrees of relations.
More detailed description
of this work as well as related materials will
be given in \S 5.
\par
On the other hand, in the case of $g=2$, 
Mess \cite{Me} proved that ${\cal I}_2={\cal K}_2$ is an infinitely
generated free group. Thus we can ask

\begin{problem} 
(i)\qua Determine whether the group $\Kg$ is finitely generated or not
for $g\geq 3$.

(ii)\qua Determine the abelianization $H_1(\Kg)$ of $\Kg$.
\label{prob:K}
\end{problem}

We mention that $\Kg$ is far from being a free group for $g\geq 3$.
This is almost clear because
it is easy to construct
subgroups of $\Kg$ which are free abelian groups of high ranks
by making use of Dehn twists along
mutually disjoint separating simple closed curves on $\Sg$.
More strongly, we can show, roughly as follows, that
the cohomological dimension of $\Kg$ will become
arbitrarily large if we take the genus $g$ sufficiently large.
Let
$$
\tau_g(2)\co \Kg\lra \hg(2)
$$
be the second Johnson homomorphism given in \cite{Mrt9, Mrt12}
(see \S 5 below for notation). 
Then it can be shown that the associated homomorphism
$$
\tau_g(2)^*\co H^*(\hg(2))\lra H^*(\Kg)
$$
is non-trivial by evaluating cohomology classes coming
from $H^*({\mathfrak h}_g(2))$, under the homomorphism
$\tau_g(2)^*$,
on abelian cycles of $\Kg$
which are supported in the above free abelian subgroups.
\par
In \S 6.6, we will consider the cohomological structure
of the group $\Kg$ from a hopefully deeper point 
of view which is related to the secondary characteristic classes
of surface bundles introduced in \cite{Mrt13}.
\par

\section{The stable cohomology of $\Mg$ and 
the stable homotopy type of ${\mathbf M}_g$}

Let $\pi\co\overline{\mathbf C}_g\ra \overline{\mathbf M}_g$ be 
the universal family of stable curves over the 
Deligne--Mumford compactification of the
moduli space ${\mathbf M}_g$. In \cite{Mu}, Mumford defined 
certain classes 
$$
\kappa_i\in A^i(\overline{\mathbf M}_g)
$$ 
in the Chow algebra (with coefficients in $\bQ$)
of the moduli space $\overline{\mathbf M}_g$
by setting $\kappa_i=\pi_*(c_1(\om)^{i+1})$ where $\om$ denotes the
relative dualizing sheaf of the morphism $\pi$.
On the other hand, in \cite{Mrt1} the author independently
defined certain integral cohomology classes 
$$
e_i\in H^{2i}(\Mg;\bZ)
$$
of the mapping class group $\Mg$
by setting $e_i=\pi_*(e^{i+1})$ where
$$
\pi\co \ED_+\Sg\ra\BD_+\Sg
$$ 
is the 
universal oriented $\Sg$--bundle and $e\in H^2(\ED_+\Sg;\bZ)$ is the
Euler class of the relative tangent bundle of $\pi$.
As was mentioned in \S 1, there exists a natural isomorphism
$H^*(\Mg;\bQ)\cong H^*({\mathbf M}_g;\bQ)$ and
it follows immediately from the definitions
that $e_i=(-1)^{i+1}\kappa_i$ as an element of these {\it rational}
cohomology groups.
The difference in signs comes from the fact that 
Mumford uses the first Chern class of the 
relative dualizing sheaf of $\pi$
while our definition
uses the Euler class of the relative tangent bundle. 
These classes $\kappa_i, e_i$ are called 
tautological classes or 
Mumford--Morita--Miller classes.
\par
In this paper, we use our notation $e_i$ to emphasize that
we consider it as an integral cohomology class of the mapping class
group rather than an element of the Chow algebra of the moduli space.
A recent work of Kawazumi and Uemura in \cite{KU} shows that
the integral class $e_i$ can play an interesting role in 
a study of certain cohomological properties of 
finite subgroups of $\Mg$.
\par
Let
\begin{equation}
\Phi\co \bQ[e_1,e_2,\cdots]\lra \lim_{g\to\infty} H^*(\Mg;\bQ)
\label{3}
\end{equation}

be the natural homomorphism from the polynomial algebra generated by
$e_i$ into the stable cohomology group of the mapping class group 
which exists by virtue of a fundamental result of Harer \cite{Har2}.
It was proved by Miller \cite{Mil} and the author \cite{Mrt2}, independently, 
that the homomorphism $\Phi$ is {\it injective} and we have the 
following well known conjecture (see Mumford \cite{Mu}).

\conjecture 
The homomorphism $\Phi$ is an isomorphism so that
$$
\lim_{g\to\infty} H^*(\Mg;\bQ)\cong \bQ[e_1,e_2,\cdots].
$$
\label{conj:basic}

We would like to mention here a few pieces of evidence which support the above
conjecture.
First of all, Harer's explicit computations in 
\cite{Har1, Har4, Har6} 
verify the conjecture in low degrees. See also
\cite{AC0} for more recent development. Secondly, Kawazumi
has shown in \cite{Kaw2} (see also \cite{Kaw1, Kaw3}) that 
the Mumford--Morita--Miller classes
occur naturally in his algebraic model of the cohomology of
the moduli space which is constructed in the framework of the complex 
analytic Gel'fand--Fuks cohomology theory, 
whereas no other classes can be obtained in this way. 
Thirdly, in \cite{KM1, KM2} Kawazumi and the author showed that
the image of the natural homomorphism
$$
H^*(H_1(\Ig);\bQ)^{Sp}\lra H^*(\Mg;\bQ)
$$
is exactly equal to the subalgebra generated by the classes $e_i$
(see \S 6.4 for more detailed survey of related works).
Here $Sp$ stands for $Sp(2g,\bZ)$.
Finally, as is explained in a survey paper by Hain and Looijenga \cite{HL} 
and also in our paper \cite{KM1},
a combination of this result with Hain's fundamental work in \cite{H3}
via Looijenga's idea to use
Pikaart's purity theorem in \cite{Pi} implies that
there are no new classes in
the continuous cohomology of $\Mg$, with respect to a certain 
natural filtration on it, in the stable range.
\par
Now there seems to be a rather canonical way of realizing 
the homomorphism $\Phi$ of
\eqref{3}
at the space level. To describe this, we first recall 
the cohomological nature of the
classical representation $\rho_0\co\Mg\ra Sp(2g,\bZ)$.
The Siegel modular group
$Sp(2g,\bZ)$ is a discrete subgroup of $Sp(2g,\bR)$ and
the maximal compact subgroup of the latter group is isomorphic
to the unitary group $U(g)$. Hence there exists a universal $g$--dimensional
complex vector bundle on the classifying space of $Sp(2g,\bZ)$.
Let $\eta$ be the pull back, under $\rho_0$, of this bundle to 
the classifying space of $\Mg$. As was explained in \cite{At1}
(see also \cite{Mrt2}), the dual
bundle $\eta^*$
can be identified, on each family $\pi\co E\ra X$
of Riemann surfaces, as follows. Namely it is the vector bundle
over the base space $X$ whose fiber on $x\in X$ is the space
of holomorphic differentials on the Riemann surface $E_x$. 
In the above paper, Atiyah used the Grothendieck Riemann--Roch theorem
to deduce the relation
$$
e_1=12 c_1(\eta^*).
$$

If we apply the above procedure to the universal family
${\mathbf C}_g\ra {\mathbf M}_g$, then we obtain a complex vector
bundle $\eta^*$ (in the orbifold sense) over ${\mathbf M}_g$ 
(in fact, more generally, over the Deligne--Mumford compactification
$\overline{\mathbf M}_g$) which is called the Hodge bundle. 
In \cite{Mu}, Mumford applied the Grothendieck Riemann--Roch theorem
to the morphism $\overline{\mathbf C}_g\ra \overline{\mathbf M}_g$
and obtained an identity, in the Chow algebra $A^*(\overline{\mathbf M}_g)$,
which expresses the  Chern classes of the Hodge bundle 
in terms of the tautological classes
$\kappa_{2i-1}$ with {\it odd} indices together with some
canonical classes coming from the boundary. 
From this identity,  we can deduce the relations
\begin{equation}
e_{2i-1}= \frac{2i}{B_{2i}}\ s_{2i-1}(\eta^*)\quad (i=1,2,\cdots)
\label{4}
\end{equation}

in the rational cohomology of $\Mg$. Here $B_{2i}$
denotes the $2i$-th Bernoulli number and $s_i(\eta^*)$ is the characteristic
class of $\eta^*$ corresponding to the formal sum $\sum_j t_j^i$ 
(sometimes called the $i$-th Newton class).
We have also obtained the above relations in \cite{Mrt1}
by applying the Atiyah--Singer index 
theorem \cite{AS} for families of elliptic operators, along the
lines of Atiyah's argument in \cite{At1}. Since $\eta^*$ is {\it flat} as 
a real vector bundle, all of its Pontrjagin classes vanish
so that we can conclude that the Chern classes of $\eta^*$ can be expressed
entirely in terms of the classes $e_{2i-1}$.
Thus we can say that the totality of the classes $e_{2i-1}$ 
of {\it odd} indices is equivalent to
the total Chern class of the Hodge bundle which comes from the Siegel
modular group.
\par
Although the rational cohomology of $\Mg$ and ${\mathbf M}_g$ 
are canonically
isomorphic to each other, there seems to be a big difference between
the torsion cohomology of them. To be more precise, let
$$
\BD_+\Sg=K(\Mg,1)\lra {\mathbf M}_g \quad (g\geq 2)
$$
be the natural mapping which is uniquely defined up to
homotopy, where
the equality above is due to a result of
Earle and Eells \cite{EE} as was already mentioned in the
introduction.
As is well known (see eg \cite{Har1}), $\Mg$ is perfect for
all $g\geq 3$ so that we can apply Quillen's plus construction
on $K(\Mg,1)$ to obtain a simply connected space $K(\Mg,1)^+$
which has the same homology as that of $\Mg$. It is known that
the moduli space ${\mathbf M}_g$ is simply connected.
Hence, by the universal property of the plus construction,
the above mapping factors through a mapping
$$
K(\Mg,1)^+\lra {\mathbf M}_g.
$$

\begin{problem}
Study the homotopy theoretical properties of 
the above mapping $K(\Mg,1)^+\lra {\mathbf M}_g$. 
In particular, 
what is its homotopy fiber ?
\end{problem}

The classical representation
$\rho_0\co\Mg\ra Sp(2g,\bZ)$ induces a mapping
\begin{equation}
K(\Mg,1)^+\lra K(Sp(2g,\bZ),1)^+
\label{5}
\end{equation}

because $Sp(2g,\bZ)$ is also perfect for $g\geq 3$.
Homotopy theoretical
properties of this map (or rather its direct limit
as $g\to\infty$) have been studied
by many authors and they 
produced interesting implications on the torsion cohomology of 
$\Mg$ (see \cite{CL1, CL2, GM, Ti} 
as well as their references). 
A final result along these lines
was obtained by Tillmann.  This says that
$K({\cal M}_{\infty},1)^+$ is an infinite loop space and
the natural map $K({\cal M}_{\infty},1)^+\ra K(Sp(2\infty,\bZ),1)^+$ 
is that of infinite loop spaces
(see \cite{Ti} for details).
See also \cite{Mrv} for a different feature of the above map,
\cite{V} for a homotopy theoretical implication of
Conjecture \ref{conj:basic} and \cite{Od} for the {\it etale} homotopy type of
the moduli spaces.
\par
Let $F_g$ be the homotopy fiber of the above mapping \eqref{5}.
Then, we have a map
$$
{\mathbf T}_g\lra F_g.
$$
Using the fact that any class $e_i$ is primitive
with respect to Miller's loop space structure on $K({\cal M}_\infty,1)^+$, 
it is easy to see that the
natural homomorphism
$$
\bQ[e_2,e_4,\cdots]\lra H^*(F_g;\bQ)
$$
is {\it injective} 
in a certain stable range
and we can ask how these cohomology classes
behave on the Torelli space.
\par
We would like to show that the classes $e_{2i}$ of {\it even}
indices are closely related to the Pontrjagin classes of
the moduli space ${\mathbf M}_g$ and also of the Torelli space
${\mathbf T}_g$. To see this, recall that ${\mathbf T}_g$ is a
complex manifold and
${\mathbf M}_g$ is {\it nearly} a complex manifold of dimension $3g-3$.
More precisely, as is well known 
it has a finite ramified covering $\widetilde{\mathbf M}_g$
which is a complex manifold and we can write
${\mathbf M}_g=\widetilde{\mathbf M}_g/G$ where $G$ is a
suitable finite group acting holomorphically on
$\widetilde{\mathbf M}_g$. Hence we have the 
Chern classes
$$
c_i\in H^{2i}(\widetilde{\mathbf M}_g;\bZ)\quad (i=1,2,\cdots)
$$
of the tangent bundle of $\widetilde{\mathbf M}_g$
which is invariant under the action of $G$.
Hence we have the rational cohomology classes  
$$
c_i^o\in H^{2i}({\mathbf M}_g;\bQ)
$$
which is easily seen to be independent of the choice of 
$\widetilde{\mathbf M}_g$.
We may call them {\it orbifold} Chern classes of the moduli space.
To identify these classes, we use the Grothendieck Riemann--Roch
theorem applied to the morphism $\pi\co {\mathbf C}_g\ra {\mathbf M}_g$
$$
\pi_*(ch(\xi)Td(\om^*))=ch(\pi_!(\xi))
$$
where $\om^*$ denotes the relative tangent bundle (in the orbifold
sense) of $\pi$ and $\xi$ is a vector bundle over ${\mathbf C}_g$.
If we take $\xi$ to be the relative cotangent bundle $\om$
as in \cite{Mu},
then we obtain the relations \eqref{4} above.
Instead of this, let us take $\xi$ to be $\om^*$. 
Since $\pi_!(\om^*)=-T{\mathbf M}_g$ by the Kodaira--Spencer
theory, we have 
\begin{align*}
ch^o({\mathbf M}_g)=&-\pi_*\bigl(ch(\om^*)Td(\om^*)\bigr)\\
=&-\pi_*\bigl(\exp\ e \ \frac{e}{1-\exp^{-e}}\bigr)\\
=& -\pi_*\Big\{\bigl(1+e+\cdots +\frac{1}{n!} e^n+\cdots\bigr)
\bigl(1+\frac{1}{2} e+\sum_{k=1}^\infty (-1)^{k-1}\frac{B_k}{(2k)!}
 e^{2k}\bigr)\Big\}
\end{align*}
where $e\in H^2({\mathbf C}_g;\bQ)$ denotes the Euler class of $\om^*$.
From this, we can conclude

\begin{align*}
s_{2k-1}^o({\mathbf M}_g)=-\Big\{&\frac{1}{(2k)!}+
\frac{1}{(2k-1)!}\cdot \frac{1}{2}
+\frac{1}{(2k-2)!}\ \frac{B_1}{2}+\frac{1}{(2k-4)!}\cdot - \frac{B_2}{4!}\\
&+\cdots +
\frac{1}{2}\cdot (-1)^k \frac{B_{k-1}}{(2k-2)!}+(-1)^{k-1}\frac{B_k}{(2k)!}
\Big\} e_{2k-1}\\
s_{2k}^o({\mathbf M}_g)=-\Big\{&\frac{1}{(2k+1)!}+
\frac{1}{(2k)!}\cdot \frac{1}{2}
+\frac{1}{(2k-1)!}\ \frac{B_1}{2}+\frac{1}{(2k-3)!}\cdot - \frac{B_2}{4!}\\
&+\cdots +
\frac{1}{6}\cdot (-1)^k \frac{B_{k-1}}{(2k-2)!}+(-1)^{k-1}\frac{B_k}{(2k)!}
\Big\} e_{2k}.
\end{align*}

The first few classes are given by
$$
s_1^o({\mathbf M}_g)=-\frac{13}{12} e_1, \quad 
s_2^o({\mathbf M}_g)=-\frac{1}{2} e_2, \quad
s_3^o({\mathbf M}_g)=-\frac{119}{720} e_3.
$$

Thus the orbifold Chern classes of ${\mathbf M}_g$ turn out to be,
in some sense, independent of $g$.
The pull back of these classes to the Torelli space ${\mathbf T}_g$ are
equal to the (genuine) Chern classes of it because ${\mathbf T}_g$ is a complex
manifold. Since the pull back of $e_{2i-1}$ to ${\mathbf T}_g$ vanishes for
all $i$, we can conclude that $s_{2i-1}({\mathbf T}_g)=0$ and
only the classes $s_{2i}({\mathbf T}_g)$ may remain to be
non-trivial. As is well known, these classes are equivalent to the
Pontrjagin classes of ${\mathbf T}_g$ as a differentiable manifold.
\par
In view of the above facts, it may be said that the classifying map
${\mathbf M}_g\ra BU(3g-3)$ 
of the holomorphic tangent bundle of
${\mathbf M}_g$ would realize the conjectural isomorphism \eqref{3}
at the space level
(rigorously speaking, we have to use some finite covering of ${\mathbf M}_g$). 
Alternatively we could use the map
$$
{\mathbf M}_g\ra {\mathbf A}_g\times BSO(6g-6)
$$ 
where the second factor
is the classifying map of the tangent bundle of ${\mathbf M}_g$
as a {\it real} vector bundle.
In short, we can say that the odd classes $e_{2i-1}$ 
serve as Chern classes of the Hodge bundle while the even classes
$e_{2i}$ embody the orbifold Pontrjagin classes of
the moduli space.
\par
According to Looijenga \cite{Lo2}, the Deligne--Mumford compactification
$\overline{\mathbf M}_g$ can also be described as a finite quotient of 
some compact complex manifold. Hence we have its orbifold Chern classes
as well as orbifold Pontrjagin classes. On the other hand,
since $\overline{\mathbf M}_g$ is a rational homology manifold, 
its combinatorial Pontrjagin classes in the sense of Thom are defined.
\par

\begin{problem}
Study the relations between orbifold Chern classes, 
orbifold Pontrjagin classes and
Thom's combinatorial Pontrjagin classes
of $\overline{\mathbf M}_g$. In particular, study the relation
between the corresponding charateristic numbers. 
\end{problem}

If we look at the basic extension \eqref{1} given in \S 2, keeping in mind
the above discussions together with the Borel vanishing theorem given in
\cite{Bo1, Bo2} concerning the triviality of twisted cohomology of
$Sp(2g,\bZ)$ with coefficients in non-trivial algebraic representations
of $Sp(2g,\bQ)$, we arrive at the following conjecture.

\conjecture
Any class $e_{2i}$ of even index is non-trivial in the 
rational cohomology of the Torelli group $\Ig$ for sifficiently large $g$.
Moreover the $Sp$--invariant part of
the rational cohomology of $\Ig$
stabilizes and we have an isomorphism
$$
\lim_{g\to\infty} H^*(\Ig;\bQ)^{Sp}\cong \bQ[e_2,e_4,\cdots].
$$
\label{conj:torellieven}

At present, even the non-triviality of the first one $e_2$ is not known.
One of the difficulties in proving this lies in the fact that 
the rational cohomology of $\Ig$ is {\it infinite} dimensional in general.
Mess observed this fact for $g=2,3$ and recently Akita \cite{Akita} proved that
$H^*(\Ig;\bQ)$ is infinite dimensional for all $g\geq 7$.  
His argument can be roughly described as follows.
He compares the orbifold Euler characteristic of ${\mathbf M}_g$
given by Harer--Zagier in \cite{HZ} with that of ${\mathbf A}_g$ given by
Harder \cite{Ha} to conclude that the Euler number of ${\mathbf T}_g$, if
defined, cannot be an integer because the latter number is much 
larger than the former one.
On the other hand, it seems to be extremely difficult to construct a family of
Riemann surfaces such that its monodromy does not act on the
homology of the fiber wheras the moduli moves in such a way that the 
classes $e_{2i}$ are non-trivial
(see a recent result of I Smith described in \cite{ABKP} for example).
Perhaps completely different approaches to this problem
along the lines of works of Jekel \cite{Je} or
Klein \cite{Kl} might also be possible.

\section{Symplectic representation theory}

As was explained in \S 2, it is an important
method of studying the structure of the mapping class group
to combine those of the Siegel modular
group $Sp(2g,\bZ)$ and the Torelli group $\Ig$ together
with the action of the former group on the structure of the latter group.
More precisely, there arise various representations of the algebraic
group $Sp(2g,\bQ)$ in the study of $\Mg$.
For example, the rational homology group $H_\bQ=H_1(\Sg;\bQ)$ 
of the surface $\Sg$ is the
fundamental representation of $Sp(2g,\bQ)$ and
Johnson's result implies that
$H_1(\Ig;\bQ)\cong \La^3 H_\bQ/H_\bQ$ 
is also a rational
representation of it. 
Hereafter, the representation $\La^3 H_\bQ/H_\bQ$ will be denoted
by $U_\bQ$.
Thus the classical
representation theory of $Sp(2g,\bQ)$ can play crucial roles.
\par
On the other hand, as was already mentioned in the introduction,
Kontsevich \cite{Ko2, Ko3} used Weyl's classical representation theory
to describe invariant tensors of various representation spaces
which appear in low dimensional topology
in terms of graphs.
In this section, we adopt this method to describe invariant tensors of
various $Sp$--modules related to the mapping class group
as well as the Torelli group.
\par
As is well known, irreducible
representations of $Sp(2g,\bQ)$ can be described
as follows (see a book by Fulton and Harris \cite{FH}). 
Let ${\mathfrak {sp}}(2g,\bC)$ be the Lie
algebra of $Sp(2g,\bC)$ and let $\mathfrak h$ be its Cartan subalgebra
consisting of diagonal matrices. Choose a system of fundamental weights
$L_i\co {\mathfrak h}\ra \bR \ (i=1,\cdots,g)$ as in \cite{FH}.
Then for each $g$--tuple
$(a_1,\cdots,a_g)$ of non-negative integers, there exists
an irreducible representation with highest weight
$(a_1+\cdots+a_g)L_1+(a_2+\cdots+a_g)L_2+\cdots+a_gL_g$. 
In \cite{FH}, this representation is denoted by $\Ga_{a_1,\cdots,a_g}$.
In this paper, following \cite{AN} we use the notation
$[a_1+\cdots +a_g, a_2+\cdots +a_{g},\cdots,a_g]$
for it. In short, irreducible
representations of $Sp(2g,\bC)$ are indexed by
Young diagrams whose number of rows are less than or equal to $g$.
These representations
are all rational representations defined over $\bQ$ so that we can consider
them as irreducible representations of $Sp(2g,\bQ)$.
For example $H_\bQ=\Ga_1=[1], \ U_\bQ=\Ga_{0,0,1}=[111]$
(which will be abbreviated by $[1^3]$ and similarly for
others with duplications)
and $S^k H_\bQ=\Ga_{k}=[k]$ where $S^k H_\bQ$ denotes the $k$-th
symmetric power of $H_\bQ$.
\par
Recall from \S 2 that
$\om_0\in H^{\otimes 2}$ denotes the symplectic class defined as
$\om_0=\sum_i (x_i\otimes y_i-y_i\otimes x_i)$ for any symplectic
basis $x_1,\cdots,x_g,y_1,\cdots, y_g$ of $H$.
As is well known, $\om_0$ is the generator of
$(H_\bQ^{\otimes 2})^{Sp}$. Also the intersection pairing
$\mu\co H\otimes H\ra \bQ$
serves as the generator of $\Hom(H_\bQ^{\otimes 2},\bQ)^{Sp}$.
\par

\subsection{Invariant tensors of $H_\bQ^{\otimes 2k}$ and its dual}

It is one of the classical results of Weyl that any 
invariant tensor of $H_\bQ^{\otimes 2k}$,
namely any element of $(H_\bQ^{\otimes 2k})^{Sp}$ can be described
as follows. A {\it linear chord diagram} $C$ with $2k$ vertices is 
a decomposition of the set of labeled vertices
$\{1,2,\cdots,2k-1,2k\}$ into pairs
$\{(i_1,j_1),(i_2,j_2),\cdots,(i_k,j_k)\}$ such that
$i_1<j_1, i_2<j_2,\cdots, i_k<j_k$ (cf Bar-Natan \cite{BN}, see also 
\cite{GN}).
We connect two vertices in each pair $(i_s,j_s)$ by an edge so that $C$
becomes a graph with $k$ edges.
We define $\sgn C$ by
$$
\sgn C=\sgn \pmatrix
1 & 2 & \cdots & 2k-1 & 2k \\
i_1 & j_1 & \cdots & i_{k} & j_{k}
\endpmatrix
.
$$
It is easy to see that there are exactly $(2k-1)!!$ linear chord
diagrams with $2k$ vertices.
For each linear chord diagram $C$, let
$$
a_C\in (H_\bQ^{\otimes 2k})^{Sp}
$$
be the invariant tensor defined by permuting the tensor product
$(\om_0)^{\otimes k}$ in such a way that the $s$-th part
$(\om_0)_s$ goes to $(H_\bQ)_{i_s}\otimes (H_\bQ)_{j_s}$, 
where $(H_\bQ)_i$ denotes
the $i$-th component of $H_\bQ^{\otimes 2k}$,
and multiplied by the factor $\sgn C$.
We also consider the dual element
$$
\al_C\in \Hom(H_\bQ^{\otimes 2k},\bQ)^{Sp}
$$
which is defined by applying the intersection pairing $\mu$
on each two components corresponding to
pairs $(i_s,j_s)$ of $C$ and multiplied by $\sgn C$. Namely we set
$$
\al_C(u_1\otimes\cdots\otimes u_{2k})=\sgn C 
\prod_{s=1}^k u_{i_s}\cdot u_{j_s}\quad
(u_i\in H_\bQ).
$$
\par
Let us write 
$$
{\cal D}^\ell(2k)=\{C_i; i=1,\cdots, (2k-1)!!\}
$$
for the set of all linear
chord diagrams with $2k$ vertices.

\begin{lemma}
$\dim (H_\bQ^{\otimes 2k})^{Sp}=\dim \Hom (H_\bQ^{\otimes 2k},\bQ)^{Sp}
=(2k-1)!!$ for $k\leq g$.
\label{lem:h2k}
\end{lemma}

\begin{proof}
Let $x_1,\cdots,x_g,y_1,\cdots,y_g$ be a symplectic basis of $H$.
There are $2g$ members in this basis while 
if $k\leq g$, then there are
only $2k \ (\leq 2g)$ positions in the tensor product $H_\bQ^{2k}$. It is
now a simple matter to construct $(2k-1)!!$ elements $\xi_j$ in  
$H_\bQ^{2k}$ such that 
$\bigl(\al_{C_i}(\xi_j)\bigr)\ (C_i\in {\cal D}^\ell(2k))$ 
is the identity matrix. Hence the elements $\{\al_{C_i}\}_i$ 
are linearly independent.
By the obvious duality, the $Sp$--invariant
components of tensors $\{a_{C_i}\}_i$ 
are also linearly independent.
\end{proof}

\begin{remark}
The stable range of the {\it $Sp$--invariant
part} of $H_\bQ^{\otimes 2k}$, which is $k\leq g$,
is twice the stable range of the irreducible decomposition of 
it, which is $k\leq\frac{g}{2}$. A similar statement is true for
other $Sp$--modules related to the mapping class group,
eg, $\La^*(\La^3 H_\bQ)$ and $\La^* U_\bQ$
(see Remark at the end of \S 4.2).
\end{remark}

Let $C, C'\in {\cal D}^\ell(2k)$ 
be two linear chord diagrams with $2k$ vertices.
Then the number $\al_C(a_{C'})$ is given by
$$
\al_C(a_{C'})=\sgn (C,C') (2g)^r
$$
where $r$ is the number of connected components of the graph
$C\cup C'$ and $\sgn (C,C')=\pm 1$ is suitably defined.
If $k\leq g$, then Lemma 4.1 above implies that
the matrix $\bigl(\al_{C_i}(a_{C_j})\bigr)$
is non-singular. 
If we go into the unstable range, degenerations occur and
it seems to be not so easy to analyze them. However, the first 
degeneration turns out to be remarkably simple and can be
described as follows.
\par

\begin{proposition}
If $g=k-1$, then the dimension of $Sp$--invariant part of $H_\bQ^{\otimes 2k}$
is exactly one less than the stable dimension. Namely
$$
\dim (H_\bQ^{\otimes 2k})^{Sp}=(2k-1)!!-1
$$
and the unique linear relation between the elements 
$a_C\ (C\in {\cal D}^\ell(2k))$
is given by
$$
\sum_{C\in{\cal D}^\ell(2k)} a_C=0.
$$
\label{prop:relation}
\end{proposition}

\begin{proof}[Sketch of proof]
For $k=1$ the assertion is empty and
for $k=2$ we can check the assertion by a direct computation.
Using the formula for the number $\al_C(a_{C'})$ given above, 
it can be shown that
$$
\sum_{C\in{\cal D}^\ell(2k)} \al_{C'}(a_C)=2^k g(g-1)\cdots (g-k+1)
$$
for any $C'\in{\cal D}^\ell(2k)$. Hence 
$\sum_{C\in{\cal D}^\ell(2k)} a_C=0$ for
$g\leq k-1$.
On the other hand, we can inductively construct $(2k-1)!!-1$
elements in $(H_\bQ^{\otimes 2k})^{Sp}$ which are linearly independent
for $g=k-1$.
\end{proof}

\begin{remark}
After we had obtained the above Proposition \ref{prop:relation}, 
a preprint by Mihailovs \cite{Mih}
appeared in which he gives a beautiful basis of $(H_\bQ^{\otimes 2k})^{Sp}$
{\it for all genera $g$}. Members of his basis are linearly ordered and the
above element 
$\sum_{C\in{\cal D}^\ell(2k)} a_C$ appears as the last one for 
$g=k$. (More precisely, his last element $\om^k$ in his notation is
equal to $k!$ times our element above.) In particular, the dimension
formula above follows immediately from his result. We expect that
we can use his basis in our approach to the Faber's conjecture
(see \S 6.4 for more details).
\end{remark}

\subsection{Invariant tensors of $\La^*(\La^3 H_\bQ)$ and $\La^* U_\bQ$}

In our paper \cite{Mrt12}, we described invariant tensors of
$\La^*(\La^3 H_\bQ)$ and $\La^* U_\bQ$ (or rather those of their
duals) in terms of trivalent graphs.
It turns out that they are specific cases of Kontsevich's
general framework given in \cite{Ko2, Ko3}.
Here we briefly summarize them. 
These descriptions were utilized in \cite{Mrt12, KM1} to construct
explicit group cocycles
for the characteristic classes $e\in H^2(\Msg;\bQ)$ and 
$e_i\in H^{2i}(\Mg;\bQ)$ (see \S 6.4 for more details). 
\par
As is well known, $\La^{2k}(\La^3 H_\bQ)$ can be considered as a 
natural quotient as well as a subspace of $H_\bQ^{\otimes 6k}$.
More precisely, let 
$p\co H_\bQ^{\otimes 6k}\ra \La^{2k}(\La^3 H_\bQ)$ be the natural
projection and let
$i\co \La^{2k}(\La^3 H_\bQ)\ra H_\bQ^{\otimes 6k}$  
be the inclusion induced from the embedding
$$
\La^3 H_\bQ\ni u_1\land u_2\land u_3\mapsto \sum_{\sigma}
\sgn \sigma\ u_{\sigma(1)}\otimes u_{\sigma(2)}\otimes u_{\sigma(3)}
\in H_\bQ^{\otimes 3} 
$$
and the similar one $\La^{2k}H_\bQ^{\otimes 3}\subset H_\bQ^{\otimes 6k}$,
where $\sigma$ runs through the symmetric group ${\mathfrak S}_3$ of degree $3$.
Then for each linear chord diagram $C\in {\cal D}^\ell(6k)$,
we have the corresponding elements
$$
p_*(a_C)\in (\La^{2k}(\La^3 H_\bQ))^{Sp},
\quad i^*(\al_C)\in \Hom(\La^{2k}(\La^3 H_\bQ),\bQ)^{Sp}.
$$

Out of each linear chord diagram $C\in {\cal D}^\ell(6k)$, 
let us construct a trivalent graph
$\vGa_C$ having $2k$ vertices as follows. We group the labeled vertices
$\{1,2,\cdots,6k\}$ of $C$ into $2k$ classes $\{1,2,3\},\{4,5,6\},\cdots,
\{6k-2,6k-1,6k\}$ and then join the three vertices belonging to
each class to a single point. This yields a trivalent graph which
we denote by $\vGa_C$. It can be easily seen that if two linear
chord diagrams $C, C'$ yield isomorphic trivalent graphs $\vGa_C, \vGa_{C'}$, 
then the corresponding elements coincide
$$
p_*(a_C)=p_*(a_{C'}),\quad i^*(\al_C)=i^*(\al_{C'}).
$$
On the other hand, it is clear that we can lift any trivalent graph
$\vGa$ with $2k$ vertices to a linear chord diagram $C$ such that
$\vGa=\vGa_C$. Hence to any such trivalent graph $\vGa$, we 
can associate invariant tensors
$$
a_\vGa\in (\La^{2k}(\La^3 H_\bQ))^{Sp}, 
\quad \al_\vGa \in \Hom(\La^{2k}(\La^3 H_\bQ),\bQ)^{Sp}
$$  
by setting $a_\vGa=p_*(a_C)$ and $\al_\vGa=\frac{1}{(2k)!}\ i^*(\al_C)$ where
$C\in {\cal D}^\ell(6k)$ is any lift of $\vGa$.

Now let ${\cal G}_{2k}$ be the set of isomorphism classes
of {\it connected} trivalent graphs with $2k$ vertices
and let ${\cal G}=\coprod_{k\geq 1} {\cal G}_{2k}$ be the
disjoint union of ${\cal G}_{2k}$ for $k\geq 1$.
Let $\bQ[a_\vGa; \vGa\in{\cal G}]$
be the polynomial algebra generated by the symbol $a_\vGa$ for 
each $\vGa\in {\cal G}$.

\begin{proposition}
The correspondence ${\cal G}_{2k}\ni \vGa\mapsto 
a_\vGa\in(\La^{2k}(\La^3 H_\bQ))^{Sp}$ defines a 
surjective algebra homomorphism
$$
\bQ[a_\vGa; \vGa\in{\cal G}]\lra (\La^*(\La^3 H_\bQ))^{Sp}
$$
which is an isomorphism in degrees $\leq \frac{2g}{3}$.
Similarly the correspondence ${\cal G}_{2k}\ni \vGa\mapsto 
\al_\vGa\in \Hom(\La^{2k}(\La^3 H_\bQ),\bQ)^{Sp}$ defines a 
surjective algebra homomorphism
$$
\bQ[\al_\vGa; \vGa\in{\cal G}]\lra \Hom(\La^*(\La^3 H_\bQ),\bQ)^{Sp}
$$
which is an isomorphism in degrees $\leq \frac{2g}{3}$. 
\qed
\label{prop:alpha}
\end{proposition}

Next we consider invariant tensors of $\La^*U_\bQ$ and its dual. 
We have a natural surjection $p\co\La^3 H_\bQ\ra U_\bQ$ and this induces
a linear map $p_*\co \La^*(\La^3 H_\bQ)\ra \La^*U_\bQ$. 
If a trivalent graph $\vGa\in {\cal G}_{2k}$ has a {\it loop}, namely an edge
whose two endpoints are the same, then clearly 
$p_*(a_\vGa)=0$. Thus let ${\cal G}^0_{2k}$ be the subset of
${\cal G}_{2k}$ consisiting of those graphs {\it without} loops
and let ${\cal G}^0=\coprod_k {\cal G}^0_{2k}$.
For each element $\vGa\in {\cal G}^0$, let $b_\vGa=p_*(a_\vGa)$.
Also let $q\co \La^3 H_\bQ\ra \La^3 H_\bQ$ be the $Sp$--equivariant 
linear map defined by $q(\xi)=\xi-\frac{1}{2g-2} C\xi\land\om_0
(\xi\in\La^3 H_\bQ)$ where $C\co\La^3 H_\bQ\ra H_\bQ$ is the contraction.
Since $q(H_\bQ)=0$, it induces a homomorphism
$q\co U_\bQ\ra \La^3 H_\bQ$ and hence 
$q\co \La^{2k}U_\bQ\ra \La^{2k}(\La^3 H_\bQ)$.
Now for each element $\vGa\in {\cal G}^0_{2k}$, let
$\beta_\vGa\co\La^{2k}U_\bQ\ra\bQ$ be defined by
$\beta_\vGa=\al_\vGa\circ q$.

\begin{proposition}
The correspondence ${\cal G}^0_{2k}\ni \vGa\mapsto 
b_\vGa\in(\La^{2k} U_\bQ)^{Sp}$ defines a 
surjective algebra homomorphism
$$
\bQ[b_\vGa; \vGa\in{\cal G}^0]\lra (\La^* U_\bQ)^{Sp}
$$
which is an isomorphism in degrees $\leq \frac{2g}{3}$.
Similarly the correspondence ${\cal G}^0_{2k}\ni \vGa\mapsto 
\beta_\vGa\in \Hom(\La^{2k}U_\bQ,\bQ)^{Sp}$ defines a 
surjective algebra homomorphism
$$
\bQ[\beta_\vGa; \vGa\in{\cal G}]\lra \Hom(\La^*U_\bQ,\bQ)^{Sp}
$$
which is an isomorphism in degrees $\leq \frac{2g}{3}$. 
\qed
\label{prop:beta}
\end{proposition}

Since $\La^3 H_\bQ\cong U_\bQ\oplus H_\bQ$, there is a natural decomposition
$$
\La^{2k}(\La^3 H_\bQ)\cong \La^{2k}U_\bQ\oplus 
(\La^{2k-1}U_\bQ\otimes H_\bQ)
\oplus\cdots\oplus (U_\bQ\otimes \La^{2k-1}H_\bQ)\oplus \La^{2k}H_\bQ
$$
and it induces that of the corresponding
$Sp$--invariant parts.
Hence we can also decompose 
the space of invariant tensors of $\La^*(\La^3 H_\bQ)$ 
and its dual according to the above splitting.
In fact, Proposition \ref{prop:beta} gives
the $\La^*U_\bQ$--part of Proposition \ref{prop:alpha}. 
We can give formulas for other parts of the above decomposition which
are described in terms of numbers of loops of trivalent graphs.
We refer to \cite{KM2} for details.  
\par

\begin{remark}
As is described in the above propositions,
the stable range of the {\it $Sp$--invariant
part} of $\La^{2k}(\La^3 H_\bQ)$ and $\La^{2k}U_\bQ$
is $2k \leq\frac{2g}{3}$. This range coincides with
Harer's improved stability range of the homology of the
mapping class group given in \cite{Har5}.
It turns out that this is far more than just an accident.
In fact, this fact will play an essential role in our approach to the
Faber's conjecture (see \S 6.4 and \cite{Mrt14} for details).
\end{remark}

\subsection{Invariant tensors of $\h1g$}

In this subsection, we fix a genus $g$ and
we write ${\cal L}_{g,1}=\oplus_k {\cal L}_{g,1}(k)$ for the free
Lie algebra generated by $H$. Also we consider
the module
$$
\h1g(k)=\Ker (H\otimes{\cal L}_{g,1}(k+1)\ra {\cal L}_{g,1}(k+2))
$$
which is the degree $k$ summand of the Lie algebra consisting of
derivations of ${\cal L}_{g,1}$ which kill
the symplectic class $\om_0\in {\cal L}_{g,1}(2)$
(see the next section \S 5 for details).
We simply write ${\cal L}^\bQ_{g,1}$ and ${\mathfrak h}^\bQ_{g,1}(k)$ for
${\cal L}_{g,1}(k)\otimes\bQ$ and $\h1g(k)\otimes\bQ$ respectively.
We show that invariant tensors of ${\mathfrak h}^\bQ_{g,1}(2k)$ or its
dual, namely
any element of ${\mathfrak h}^\bQ_{g,1}(2k)^{Sp}$ or 
$\Hom({\mathfrak h}^\bQ_{g,1}(2k),\bQ)^{Sp}$
can be represented by a linear combination of chord diagrams
with $(2k+2)$ vertices. Here a {\it chord diagram} with $2k$ vertices
is a partition of $2k$ vertices lying on a circle
into $k$ pairs where each pair is connected by a chord.  
Chord diagrams already appeared in the theory of Vassiliev knot 
invariants (see \cite{BN}) and they played an important role. 
In the following, we will see that
they can play another important role also in our theory.
\par
To show this,
we recall a well known characterization of elements of 
${\cal L}_{g,1}^\bQ(k)$ in 
$H_\bQ^{\otimes k}$. There are several such characterizations
which are given in terms of various projections 
$H_\bC^{\otimes k}\ra {\cal L}_{g,1}(k)\otimes\bC$
(see \cite{Re}). Here we adopt the following one.

\begin{lemma}
Let ${\mathfrak S}_k$ be the symmetric group of degree
$k$ and let $\sigma_i=(12\cdots i)$
$\in{\mathfrak S}_k$ be the cyclic permutation.
Let $p_k=(1-\sigma_k)(1-\sigma_{k-1})\cdots (1-\sigma_2)\in\bZ[{\mathfrak S}_k]$
which acts linearly on $H_\bQ^{\otimes k}$.
Then $p_k^2=k p_k$ and an element $\xi\in H_\bQ^{\otimes k}$ belongs to
${\cal L}_{g,1}^\bQ(k)$ if and only if $p_k(\xi)=k \xi$. Moreover 
${\cal L}_{g,1}^\bQ(k)=\im p_k$.
\qed
\label{lem:ell}
\end{lemma}

If we consider ${\cal L}_{g,1}^\bQ(k+1)$ as
a subspace of $H_\bQ^{\otimes (k+1)}$, then the bracket operation
$$
H_\bQ\otimes {\cal L}_{g,1}^\bQ(k+1)\lra {\cal L}_{g,1}^\bQ(k+2)
$$
is simply given by the correspondence
$u\otimes\xi\mapsto u\otimes\xi-\xi\otimes u
\ (u\in H_\bQ, \xi\in{\cal L}_{g,1}^\bQ(k+1))$. 
Hence it is easy to deduce the following characterization of 
${\mathfrak h}^\bQ_{g,1}(k)$ inside $H_\bQ^{\otimes (k+2)}$.

\begin{proposition}
An element $\xi\in H_\bQ^{\otimes (k+2)}$ belongs to 
${\mathfrak h}^\bQ_{g,1}(k)\subset H_\bQ^{\otimes (k+2)}$
if and only if the following two conditions are satisfied.
(i) $(1\otimes p_{k+1}) \xi=(k+1) \xi$ and
(ii) $\sigma_{k+2} \xi=\xi$.
\qed
\label{prop:hg}
\end{proposition}

We can construct
a basis
of $(H_\bQ\otimes {\cal L}_{g,1}^\bQ(2k+1))^{Sp}$ as follows.
Recall that we write ${\cal D}^\ell(2k)$ 
for the set of linear chord diagrams with $2k$ vertices
so that it gives a basis of $(H_\bQ^{\otimes 2k})^{Sp}$ for $k\leq g$
(see Lemma \ref{lem:h2k}). By Lemma \ref{lem:ell}, we have
$$
H_\bQ\otimes {\cal L}_{g,1}^\bQ(2k+1)=\im (1\otimes p_{2k+1})
$$
where we consider $1\otimes p_{2k+1}$ as an endomorphism of
$H_\bQ\otimes (H_\bQ\otimes H_\bQ^{\otimes 2k})$.
Let $C_0$ be the edge which connects the first two  
of the $(2k+2)$ vertices corresponding to 
$H_\bQ\otimes (H_\bQ\otimes H_\bQ^{\otimes 2k})$.
For each element $C\in{\cal D}^\ell(2k)$, consider the disjoint union
$\widetilde C=C_0\coprod C$ which is a
linear chord diagram with $(2k+2)$ vertices.
Hence we have the corresponding invariant tensor
$$
a_{\widetilde C}\in (H_\bQ\otimes (H_\bQ\otimes H_\bQ^{\otimes 2k}))^{Sp}.
$$
Let $\ell_C=1\otimes p_{2k+1}(a_{\widetilde C})$. Then by 
Proposition \ref{prop:hg},
$\ell_C$ is an element of $(H_\bQ\otimes {\cal L}_{g,1}^\bQ(2k+1))^{Sp}$.

\begin{proposition}
If $k\leq g$, then the set of elements $\{\ell_C; C\in{\cal D}^\ell(2k)\}$
forms a basis of invariant tensors of
$H_\bQ\otimes {\cal L}_{g,1}^\bQ(2k+1)$. In particular
$$
\dim (H_\bQ\otimes {\cal L}_{g,1}^\bQ(2k+1))^{Sp}=(2k-1)!!
$$
\label{prop:ec}
\end{proposition}

\begin{proof}[Sketch of Proof]
It can be shown that the elements in $\{\ell_C; C\in{\cal D}^\ell(2k)\}$ are
linearly independent because if we express $\ell_C$ as a linear
combination of the standard basis of $(H_\bQ^{\otimes (2k+2)})^{Sp}$
given in Lemma \ref{lem:h2k},
then we find
$\ell_C=a_{\widetilde C}+\text{other terms}$.
On the other hand, we can show that the projection under
$1\otimes p_{2k+1}$ of any member of this standard basis
can be expressed as a linear combination of $\ell_C$.
\end{proof}
\par

Now we consider invariant tensors of ${\mathfrak h}^\bQ_{g,1}(2k)$.
Associated to any element $C\in {\cal D}^\ell(2k)$ we
have the corresponding 
invariant tensor $\ell_C\in (H_\bQ\otimes {\cal L}_{g,1}^\bQ(2k+1))^{Sp}$.
Consider the element
$$
\xi_C= \sum_{i=1}^{2k+2} \sigma_{2k+2}^i \ell_C.
$$
By Proposition \ref{prop:hg} (in particular condition (ii)), $\xi_C$
belongs to ${\mathfrak h}^\bQ_{g,1}(2k)$ and it is clear from the above
argument that these elements span the whole invariant space
${\mathfrak h}^\bQ_{g,1}(2k)^{Sp}$. More precisely the cyclic group
$\bZ/(2k+2)$ of order $2k+2$ acts naturally on
$(H_\bQ\otimes {\cal L}_{g,1}^\bQ(2k+1))^{Sp}$
and ${\mathfrak h}^\bQ_{g,1}(2k)^{Sp}$ is nothing but the invariant
subspace of this action.
Although there does not seem to
exist any simple formula for the dimension of this invariant
subspace, this procedure gives a method of enumerating
all the elements of it.
\par
Here is another approach to this problem
which might be practically better than
the above, in particular in the dual setting. 
We simply use two conditions (i), (ii) in Proposition \ref{prop:hg}
in the opposite way. Namely we first consider condition (ii).
Recall that any linear chord diagram $C$ with $(2k+2)$ vertices
gives rise to an $Sp$--invariant map
$$
\al_C\co H_\bQ^{\otimes (2k+2)}\lra \bQ.
$$
Let us write $C_i$ for $\sigma_{2k+2}^i C\ (i=1,\cdots ,2k+2)$.
Then, in view of condition (ii) above, the restrictions of 
$\al_{C_i}$ to the subspace 
${\mathfrak h}^\bQ_{g,1}(2k)\subset H_\bQ^{\otimes (2k+2)}$
are equal to each other for all $i$. This means that, instead of
linear chord diagram $C$, we may assume that all of the vertices
of $C$ are arranged on a circle. 
But then we obtain a usual chord diagram.
Thus we can say that any chord diagram $C$ with
$(2k+2)$ vertices defines an
element of $\Hom({\mathfrak h}^\bQ_{g,1}(2k),\bQ)^{Sp}$.
In the case of Vassiliev knot invariant, the linear space 
spanned by chord diagrams with $2k$ vertices,
which is denoted by ${\cal G}_k{\cal D}^c$ in \cite{BN}, modulo
the (4T) relation serves as the set of Vassiliev invariants
of order $k$. In our case, the linear space 
spanned by chord diagrams
with $(2k+2)$ vertices, namely ${\cal G}_{k+1}{\cal D}^c$,
modulo
the relations coming from condition (i) above can be identified
with $\Hom({\mathfrak h}^\bQ_{g,1}(2k),\bQ)^{Sp}$ (and also its dual).
In particular, we have a surjection
$$
{\cal G}_{k+1}{\cal D}^c\lra\Hom({\mathfrak h}^\bQ_{g,1}(2k),\bQ)^{Sp}.
$$
\par
In \S 6.2, we will show experimental results 
which have been obtained by explicit computations 
applying this method.
\par

\section{Graded Lie algebras related to the mapping class group}

In this section, we introduce various graded Lie algebras which are
related to the mapping class group.
We begin by recalling the definition of the Johnson homomorphisms
briefly (see \cite{Jo4, Mrt9, Mrt12} for details).
Let $\M1g$ be the mapping class group of $\Sg$ relative
to an embedded disk $D^2\subset \Sg$ 
and let $\Msg$ be the mapping class group of $\Sg$ 
relative to the base point
$*\in D^2$. We write $\Sigma_g^0$ for $\Sg\setminus\Int D^2$
so that $\pi_1\Sigma_g^0$ is a free group of rank $2g$.
Let $\zeta\in\pi_1\Sigma_g^0$ be the element represented by
a simple closed curve on $\Sigma_g^0$
which is parallel to the boundary.
Then as is well known, we have natural isomorphisms
due originally to Nielsen
\begin{align*}
\Mg&\cong\Out_+\pi_1\Sg,\quad \Msg\cong \Aut_+\pi_1\Sg\\
&\M1g\cong \{\varph\in \Aut \pi_1 \Sigma_g^0; \varph(\zeta)=\zeta\}.
\end{align*}

By virtue of this, any filtration on the fundamental groups 
of surfaces induces those of the corresponding mapping class groups.
In particular, the lower central series induces natural
filtrations. More precisely, let $\vGa_k(G)$ denote the $k$-th
term in the lower central series of a group $G$, where
$\vGa_0(G)=G$ and $\vGa_k(G)=[G,\vGa_{k-1}(G)]$ for $k\geq 1$.
We set
\begin{align*}
\M1g(k)&=\{\varph\in \M1g; \varph(\ga)\ga\inv\in\vGa_k(\pi_1\Sigma_g^0)
\ \text{for any}\ \ga\in\pi_1\Sigma_g^0\}\\
\Msg(k)&=\{\varph\in \Msg; \varph(\ga)\ga\inv\in\vGa_k(\pi_1\Sg)
\ \text{for any}\ \ga\in\pi_1\Sg\}
\end{align*}
and
$$
\Mg(k)=\pi(\Msg(k))
$$
where $\pi\co\Msg\ra\Mg$ is the natural projection. 
Thus we obtain a natural filtration $\{{\cal M}(k)\}_k$ on each of the 
three types of mapping class groups. It is easy to
see that the first one ${\cal M}(1)$ is nothing but the
Torelli group, namely $\I1g, \Isg$ or $\Ig$.
Now we can say that the Johnson homomorphism 
is the one which describes the associated graded quotients
of the mapping class groups, with respect to these fitrations, 
explicitly in terms of derivations of graded Lie algebras
associated to the lower central series of the fundamental
groups of surfaces.
\par
Let us simply write $H$ for $H_1(\Sg;\bZ)$ as before and let
${\cal L}_{g,1}=\oplus_k {\cal L}_{g,1}(k)$ be the 
free graded Lie algebra generated
by $H$. As is well known, ${\cal L}_{g,1}$ is the graded Lie algebra
associated to the lower central series of $\pi_1\Sigma_g^0$,
namely we have natural isomorphisms
${\cal L}_{g,1}(k)\cong \vGa_{k-1}(\pi_1\Sigma_g^0)/\vGa_{k}(\pi_1\Sigma_g^0)$
(see \cite{MKS}).
Let $\om_0\in \La^2 H={\cal L}_{g,1}(2)$
be the symplectic class and let $I=\oplus_{k\geq 2} I_k$ be the
ideal of ${\cal L}_{g,1}$ generated by $\om_0$.
Then a result of Labute \cite{La} says that
the quotient Lie algebra
${\cal L}_g={\cal L}_{g,1}/I$ serves as 
the graded Lie algebra
associated to the lower central series of $\pi_1\Sg$.
Now we define
\begin{align*}
\h1g(k)=& \Ker (H\otimes{\cal L}_{g,1}(k+1)\ra{\cal L}_{g,1}(k+2))\\
\cong& \Ker (\Hom(H,{\cal L}_{g,1}(k+1))\ra {\cal L}_{g,1}(k+2))\\
\hsg(k)=& \Ker (H\otimes{\cal L}_g(k+1)\ra{\cal L}_g(k+2))\\
\cong& \Ker (\Hom(H,{\cal L}_g(k+1))\ra {\cal L}_g(k+2))
\end{align*}
where the second isomorphisms in each of the terms above are induced
by the Poincar\'e duality $H^*\cong H$.
In our previous papers \cite{Mrt9, Mrt12}, ${\cal L}_{g,1}, {\cal L}_g$
have been denoted by ${\cal L}^0, {\cal L}$ and also
$\h1g(k), \hsg(k)$
have been denoted by ${\cal H}_k^0, {\cal H}_k$, respectively.

Then the $k$-th Johnson homomorphisms 
\begin{align*}
\tau_{g,1}(k)\co & \M1g(k)\lra \h1g(k)\\
\tau_{g,*}(k)\co & \Msg(k)\lra \hsg(k)
\end{align*}
are defined by the correspondence
$$
\M1g(k)\ni \varph\mapsto \tau_{g,1}(k)(\varph)=
\{[\ga]\mapsto [\varph(\ga)\ga\inv]\}
\in \Hom(H,{\cal L}_{g,1}(k+1))
$$
$(\ga\in \pi_1\Sigma_g^0)$
for $\M1g$ and similarly for $\Msg$. Here $[\ga]\in H$ denotes the
homology class of $\ga\in \pi_1\Sigma_g^0$ and $[\varph(\ga)\ga\inv]$
denotes the class of $\varph(\ga)\ga\inv\in \vGa_{k}(\pi_1\Sigma_g^0)$
in ${\cal L}_{g,1}(k+1)=\vGa_{k}(\pi_1\Sigma_g^0)/\vGa_{k+1}(\pi_1\Sigma_g^0)$.
It can be shown that 
$$
\Ker \tau_{g,1}(k)=\M1g(k+1),\quad \Ker\tau_{g,*}(k)=\Msg(k+1)
$$
so that we have isomorphisms
\begin{align*}
\M1g(k)/\M1g(k+1)\cong& \im\tau_{g,1}(k)\subset{\mathfrak h}_{g,1}(k)\\
\Msg(k)/\Msg(k+1)\cong& \im\tau_{g,*}(k)\subset{\mathfrak h}_{g,*}(k).
\end{align*}

Next we consider the filtration $\{\Mg(k)\}_k$ 
of the usual mapping class group $\Mg$.
As is mentioned above, $\Mg(k)=\pi(\Msg(k))$ and it was proved
in \cite{AK} that $\pi_1\Sg\cap\Msg(k)=\vGa_{k-1}(\pi_1\Sg)$.
Hence we have a natural injection ${\cal L}_g(k)\subset {\mathfrak h}_{g,*}(k)$.
We define the $k$-th Johnson homomorphism
$$
\tau_g(k)\co \Mg(k)\lra \hg(k)
$$
by setting
$\tau_g(k)(\varph)=\tau_{g,*}(k)(\tilde\varph)
\mod {{\cal L}_g}(k)\ 
(\varph\in\Mg(k))$
where $\hg(k)=\hsg(k)/{\cal L}_{g}(k)$ and
$\tilde\varph\in\Msg(k)$ is any lift of $\varph$.
\par
Thus we obtain three graded modules
$$
\h1g=\bigoplus_{k=1}^\infty \h1g(k),\quad
\hsg=\bigoplus_{k=1}^\infty \hsg(k),\quad
\hg=\bigoplus_{k=1}^\infty \hg(k)
$$
and it turns out that they have natural structures of graded Lie
algebras over $\bZ$. 
The relations between these three graded Lie algebras
$\h1g, \hsg,\hg$ are described simply by the following
two short exact sequences
\begin{align*}
0\lra {\mathfrak j}_{g,1}\lra &\h1g\lra\hsg\lra 0\\
0\lra {\cal L}_g\lra&\hsg\lra\hg\lra 0
\end{align*}
where 
$$
{\mathfrak j}_{g,1}=\bigoplus_{k=2}^\infty{\mathfrak j}_{g,1}(k),\quad 
{\mathfrak j}_{g,1}(k)=\Ker (\Hom(H,I_{k+1})\ra I_{k+2}))
$$
and ${\cal L}_g={\cal L}_{g,1}/I$ is the graded Lie algebra associated to the
lower central series of $\pi_1\Sg$ as before.
\par
Sometimes it is useful to consider the tensor products  with $\bQ$ of the
modules appearing above. Let us denote them by
attaching a superscript $\bQ$ to the original $\bZ$--form.
For example ${\cal L}_{g,1}^\bQ={\cal L}_{g,1}\otimes\bQ$ and
${\mathfrak h}^\bQ_{g,1}=\h1g\otimes\bQ$ which we already used in \S 4.
In these terminologies, ${\mathfrak h}^\bQ_{g,1}, {\mathfrak h}^\bQ_{g,*}$
are nothing but the graded Lie algebras
consisting of derivations of ${\cal L}_{g,1}^\bQ, {\cal L}_g^\bQ$
with positive degrees which kill the symplectic class
$\om_0$ and ${\mathfrak h}^\bQ_g$ is equal to the quotient of 
${\mathfrak h}^\bQ_{g,*}$ by
inner derivations.
We omit the degree 0 part
${\mathfrak h}^\bQ_{g,1}(0)={\mathfrak h}^\bQ_{g,*}(0)={\mathfrak h}^\bQ_g(0)
={\mathfrak{sp}}(2g,\bQ)$ because
it is the Lie algebra of the rational form of
$\M1g/\M1g(1)=\Msg/\Msg(1)=\Mg/\Mg(1)=Sp(2g,\bZ)$.  
\par
Now consider projective limits of 
nilpotent groups
\begin{align*}
{\cal N}_{g,1}&=\varprojlim_{k\geq 1} \I1g/\M1g(k),\quad 
{\cal N}_{g,*}=\varprojlim_{k\geq 1} \Isg/\Msg(k)\\
{\cal N}_g&=\varprojlim_{k\geq 1} \Ig/\Mg(k)
\end{align*}
which are associated to the filtrations ${\cal M}(k)$ on the
corresponding mapping class groups.
We can tensor these groups with $\bQ$ to obtain pronilpotent
Lie groups ${\cal N}_{g,1}^\bQ, {\cal N}_{g,*}^\bQ,
{\cal N}_g^\bQ$. Let us write
${\mathfrak n}_{g,1}, {\mathfrak n}_{g,*}, {\mathfrak n}_g$
for their Lie algebras and also let 
$\Gr{\mathfrak n}_{g,1}, \Gr{\mathfrak n}_{g,*},
\Gr{\mathfrak n}_g$ be their associated graded Lie algebras,
respectively.
Then the Johnson homomorphism induces embeddings of graded Lie algebras
$$
\Gr{\mathfrak n}_{g,1}\subset{\mathfrak h}^\bQ_{g,1},\quad 
\Gr{\mathfrak n}_{g,*}\subset{\mathfrak h}^\bQ_{g,*},
\quad
\Gr{\mathfrak n}_g\subset{\mathfrak h}^\bQ_g
$$
(in fact, these embeddings are defined at the level of $\bZ$--forms).
More precisely, we can identify ${\mathfrak n}(k)$ with $\im\tau^\bQ(k)$
so that we can write $\Gr{\mathfrak n}=\im\tau^\bQ\subset{\mathfrak h}^\bQ$
for any type of decorations $\{g,1\},\{g,*\},g$.
It is a very important problem to idetify these Lie subalgebras
inside the Lie algebras of derivations.
\par
The cohomological structure of the first type ${\mathfrak h}^\bQ_{g,1}$ 
of the above three 
graded Lie algebras was investigated by Kontsevich
in his celebrated papers \cite{Ko2, Ko3} 
where he considered three types of Lie algebras
$\ell_g, a_g, c_g$ which consist of derivations
of certain Lie, associative,
and commutative algebras. In fact our ${\mathfrak h}^\bQ_{g,1}$ is nothing but
the Lie subalgebra
$\ell_g^+$ of $\ell_g$
consisting of elements with positive degrees.
There are natural injections
$\h1g\subset{\mathfrak h}_{g+1,1}$
so that we can make the direct limit
$$
{\mathfrak h}_\infty^\bQ=\lim_{g\to\infty} {\mathfrak h}^\bQ_{g,1}
$$
which is equal to the positive part $\ell^+_\infty$
of $\ell_\infty$ in
Kontsevich's notation.
In \S 6.5, we will apply one of Kontsevich's results in the above cited
papers to our ${\mathfrak h}_\infty^\bQ$ and obtain definitions of certain
(co)homology classes of outer automorphism groups $\Out F_n$ of 
free groups $F_n$ of rank $n\geq 2$.
\par
The second graded Lie algebra, which we consider in this paper, is the
Torelli Lie algebra which is, by definition, the Malcev Lie algebra
of the Torelli group. The structure of the Torelli Lie algebra has been
extensively studied by Hain in \cite{H1, H3}. Here we summarize his results
briefly for later use in \S 6 (see the above papers for details).
We write ${\mathfrak t}_{g,1}, {\mathfrak t}_{g,*}, {\mathfrak t}_g$
for the Torelli Lie algebras which correspond to three types of the 
Torelli groups $\I1g, \Isg,\Ig$, respectively (Hain uses the notation
${\mathfrak t}^1_g$ for ${\mathfrak t}_{g,*}$). 
\par
In the above, we considered certain surjective homomorphisms
$$
\I1g\lra \N1g,\quad \Isg\lra \Nsg,\quad \Ig\lra \Ng
$$ 
from each type of the Torelli groups to a tower of 
torsion free nilpotent groups. 
Roughly speaking, the Malcev completion of the Torelli group 
(or more generally of any finitely generated group) is defined to be the
projective limit of such homomorphisms. Since any finitely generated
torsion free nilpotent group $N$ can be canonically embedded into its 
Malcev completion $N\otimes\bQ$ which is a Lie group over $\bQ$,
the Malcev completion of the Torelli groups can be described by
certain homomorphisms
$$
\I1g\lra \T1g,\quad \Isg\lra \Tsg,\quad \Ig\lra \Tg
$$ 
from the Torelli groups into pronilpotent Lie groups over $\bQ$. They
are characterized by the universal property that for any homomorphism
$\rho\co {\cal I}\ra {\cal N}$ from the Torelli group into 
a pronilpotent group ${\cal N}$, 
there exists a unique homomorphism $\bar\rho\co {\cal T}\ra {\cal N}$ 
such that $\rho=\bar\rho\circ\mu$ where $\mu\co {\cal I}\ra{\cal T}$ is the
homomorphism given above (we omit the subscripts).
Since any nilpotent Lie group is determined by its Lie algebra,
the pronilpotent group ${\cal T}$ is determined by its Lie algebra
${\mathfrak t}$ which is a pronilpotent Lie algebra over $\bQ$.
These are the definitions of the pronilpotent Lie algebras
${\mathfrak t}_{g,1},\tsg,\tg$ which we would like to call the Torelli Lie 
algebras.
Let $\Gr {\mathfrak t}=\oplus_k {\mathfrak t}(k)$ be the graded Lie algebra 
associated 
to the lower central series of ${\mathfrak t}$ and also let $\Gr {\cal I}$ be
the graded Lie algebra (over $\bZ$) associated to the lower central series
of the Torelli group ${\cal I}$. Then a general fact
about the Malcev completion implies that 
there is a natural isomorphism $\Gr {\mathfrak t}\cong (\Gr {\cal 
I})\otimes\bQ$.
In particular, we have an isomorphism
$$
{\mathfrak t}(k)\cong (\vGa_{k-1}({\cal I})/\vGa_k({\cal I})) \otimes\bQ.
$$
\par
Now by the universal property of the Malcev completion, there is a
uniquely defined homomorphism ${\cal T}\ra {\cal N}$ which induces
a morphism ${\mathfrak t}\ra {\mathfrak n}$. This induces a homomorphism
$\Gr{\mathfrak t}\ra \Gr{\mathfrak n}$ of the associated graded Lie algebras. 
Thus we obtain homomorphisms
\begin{align*}
{\mathfrak t}_{g,1}(k)\ra \n1g(k)&\subset {\mathfrak h}^\bQ_{g,1}(k),\quad 
{\mathfrak t}_{g,*}(k)\ra \nsg(k)\subset{\mathfrak h}^\bQ_{g,*}(k),\\
&{\mathfrak t}_g(k)\ra\ng(k)\subset{\mathfrak h}^\bQ_g(k).
\end{align*}

In fact, Johnson already observed in \cite{Jo4} that the $(k-1)$-th term
$\vGa_{k-1}(\I1g)$ of the lower central series of $\I1g$ is contained
in $\M1g(k)$ so that there is a natural homomorphism
$$
\vGa_{k-1}(\I1g)/\vGa_k(\I1g)\lra \M1g(k)/\M1g(k+1)\subset\h1g(k)
$$
which is defined over $\bZ$. 
\par
The third (and the final) Lie algebra, denoted by ${\mathfrak u}_g$,
is the one introduced by Hain in \cite{H1}. It lies, in some sense,
between the Torelli Lie algebra ${\mathfrak t}$ and ${\mathfrak n}$.
We have the following commutative diagram
\begin{equation}
\begin{CD}
1 @>>> \Ig @>>> \Mg @>>> Sp(2g,\bZ) @>>> 1 \\
@. @VVV@VVV@VVV \\
1 @>>> {\cal N}_g\otimes\bQ @>>> {\cal G}_g @>>> Sp(2g,\bQ) @>>> 1  \\
\end{CD}
\label{6}
\end{equation}
where ${\cal G}_g$ is the 
Zariski closure of the image of $\varprojlim \Mg/\Mg(k)$
in the automorphism group of the Malcev Lie algebra
of $\pi_1\Sg$.
Hain applied the relative Malcev completion (see \cite{H4}),
which is due to Deligne, to the mapping class group $\Mg$
relative to the classical representation
$\Mg\ra Sp(2g,\bQ)$. Roughly speaking, it is the projective
limit of representations like in \eqref{6} where ${\cal G}_g$
is replaced by an algebraic group defined over $\bQ$ such that
the image of $\Mg$ there is Zariski dense and 
$N_g\otimes\bQ$
is replaced by a unipotent subgroup.
It has the form
$$
1\lra {\cal U}_g\lra {\cal E}_g\lra Sp(2g,\bQ)\lra 1
$$
where ${\cal E}_g$ is a proalgebraic group and ${\cal U}_g$
is a prounipotent group. ${\cal U}_g$ is
called the prounipotent radical
of the relative Malcev completion of the mapping class group.
Let ${\mathfrak u}_g$ be the Lie algebra of ${\cal U}_g$.
There are also decorated versions of these groups. 
By the definition, we have canonical homomorphisms
$$
{\cal T}_g\lra {\cal U}_g\lra {\cal N}\otimes \bQ
$$
and their Lie algebra version
$$
{\mathfrak t}_g\lra {\mathfrak u}_g\lra {\mathfrak n}_g.
$$
\par
Now we can state two fundamental results of Hain as follows.
The first one gives a complete description
of the relation between the Lie algebras ${\mathfrak t}$ and ${\mathfrak u}$.

\begin{theorem}[Hain \cite{H1}]
Assume that  $g\geq 3$. Then the natural homomorphisms 
${\cal T}_{g,1}\ra {\cal U}_{g,1}$, ${\cal T}_{g,*}\ra {\cal U}_{g,*}$,
${\cal T}_g\ra {\cal U}_g$ are all surjective
and the kernel of each is a central subgroup isomorphic to $\bQ$.
\qed
\end{theorem}

Equivalent statement can be given for the associated Lie algebras
${\mathfrak t}, {\mathfrak u}$.
The second result is the explicit
presentation of the Torelli Lie algebra ${\mathfrak t}$.
Let $U=\La^3 H/H$ and let $U_\bQ$ be $U\otimes\bQ$
as before.
Then the Johnson homomorphism
$$
\tau_g^\bQ(1)\co \Ig\lra U_\bQ
$$
induces a homomorphism between the second cohomology groups
$$
\tau_g^\bQ(1)^*\co \La^2 U_\bQ\lra H^2(\Ig;\bQ)
$$
and Hain determined the kernel of this homomorphism as follows.

\begin{proposition}[Hain \cite{H3}]
Let $\tau_g^\bQ(1)\co \Ig\lra U_\bQ$ be as above. Then
$$
\Ker\ \tau_g^\bQ(1)^*=[2^2]+[0]\subset \La^2U_\bQ=
[1^6]+[1^4]+[1^2]+[2^21^2]+[2^2]+[0]
$$
where the decomposition of $\La^2U_\bQ$ is valid for $g\geq 6$.
\qed
\end{proposition}
\par
The fact that $[2^2]+[0]$ is contained in the kernel was essentially proved
in our papers \cite{Mrt5, Mrt6}, though in a more primitive way. 
In fact, we obtained the secondary invariant
$d_1\in H^1(\Kg;\bQ)$ (see \S 6.6 below for details) from
the vanishing of the trivial summand $[0]$ in $H^2(\Ig;\bQ)$
which turned out to be closely related to the Casson invariant.
Hain proved the above result in a systematic
way. In particular, the non-triviality 
in $H^2(\Ig;Q)$ of other summands than $[2^2]+[0]$
was shown by
constructing explicit abelian cycles of $\Ig$
on which they take non-zero values.
\par

\begin{theorem}[Hain \cite{H3}]
If $g\geq 3$, then ${\mathfrak t}_g$ is isomorphic to the completion of its
associated graded Lie algebra $\Gr {\mathfrak t}_g$. Moreover if $g\geq 6$, then
it has a presentation
$$
\Gr {\mathfrak t}_g={\cal L}(U_\bQ)/([1^6]+[1^4]+[1^2]+[2^21^2])
$$
where ${\cal L}(U_\bQ)$ denotes the free Lie algebra generated by $U_\bQ=[1^3]$.
\qed
\end{theorem}
\par 
Hain obtained similar presentations for other Lie algebras 
${\mathfrak t}_{g,1}, {\mathfrak t}_{g,*}, {\mathfrak u}_{g,1},{\mathfrak 
u}_{g,*}, {\mathfrak u}_g$
as well.
Although these presentations are very simple, the proof requires 
deep and powerful techniques in Hodge theory.
One of the main ingredients is to put mixed Hodge structures on ${\mathfrak 
u}_g$
and then on ${\mathfrak t}_g$ which depend on a fixed complex structure
on the reference surface. One important consequence of this is that
after tensoring with $\bC$, they are canonically isomorphic to their
associated graded Lie algebras.
\par
Thus for the study of the structure of the mapping class group,
it is enough to investigate the morphisms
$$
{\mathfrak t}_g(k)\lra {\mathfrak u}_g(k)\lra \im\tau_g^\bQ(k)
$$
between the three graded Lie algebras 
$\Gr {\mathfrak t}_g, \Gr{\mathfrak u}_g$ and
$\oplus_k\im\tau_g^\bQ(k)\subset {\mathfrak h}^\bQ_g$. 
Another important consequence of Hain's work is that 
$\im \tau_g^\bQ=\oplus_k\im\tau_g^\bQ(k)$ is generated
by degree one summand $\im\tau_g^\bQ(1)=[1^3]$.
\par

\section{A prospect on the structure of the mapping class group}
\label{prospect}

In this section, we would like to describe some of the
various aspects of the structure of the mapping class group
which seem to deserve further investigation in the future.
More precisely, we consider the following topics.
\ppar

6.1 $\;$ Structure of the graded Lie algebra $\h1g$
\par
6.2 $\;$ Description of the $Sp$--invariant part of ${\mathfrak h}^\bQ_{g,1}(6)$
\par
6.3 $\;$ Kernel of ${\mathfrak t}_g\ra\im\tau^\bQ_g$ and invariants of 
$3$--manifolds
\par
6.4 $\;$ Cocycles for the Mumford--Morita--Miller classes
\par
6.5 $\;$ Cohomology of the graded Lie algebra ${\mathfrak h}^\bQ_\infty$ and 
$\Out F_n$
\par
6.6 $\;$ Secondary characteristic classes of surface bundles
\par
6.7 $\;$ Bounded cohomology of the mapping class group
\par
6.8 $\;$ Representations of the mapping class group

\subsection{Structure of the graded Lie algebra $\h1g$}

We consider the structure of the graded Lie algebra
$$
{\mathfrak h}_{g,1}=\bigoplus_{k=1}^\infty {\mathfrak h}_{g,1}(k)
$$
(see \S 5). Recall that this is the graded Lie algebra
consisting of derivations, with positive degrees,
of the free graded Lie algebra
${\cal L}_{g,1}=\oplus_k {\cal L}_{g,1}(k)$ which kill the symplectic class 
$\om_0\in {\cal L}_{g,1}(2) =\La^2 H$.
We omit the degree 0 part $\h1g(0)\otimes\bQ$ which is 
isomorphic to ${\mathfrak sp}(2g,\bQ)$.
\par
As was explained in \S 5, this Lie algebra 
serves as the target of the Johnson homomorphisms
\begin{equation*}
\tau_{g,1}(k)\co \M1g(k)/\M1g(k+1)\lra {\mathfrak h}_{g,1}(k)
\end{equation*}
so that the main problem is to describe the image 
\begin{equation*}
\im \tau_{g,1}\subset {\mathfrak h}_{g,1}
\end{equation*}
as a Lie subalgebra of $\h1g$.
It is one of the basic results of Johnson that
\begin{equation*}
\im\tau_{g,1}(1)=\h1g(1)=\La^3 H.
\end{equation*}
Here we mention two important topological applications  
obtained by Hain in his fundamental work in \cite{H3}.
One is that the graded Lie algebra $\im\tau_{g,1}^\bQ$ 
(the superscript $\bQ$ will mean that
we take the operation
$\otimes\bQ$ as before)
is generated by the degree one summand 
$\im\tau_{g,1}^\bQ(1)=\La^3 H_\bQ$, which was already mentioned in \S 5. 
The other is that the natural map
$\M1g(k)/\M1g(k+1)\ra \Msg(k)/\Msg(k+1)$ is surjective after tensoring
with $\bQ$ (and hence is an isomorphism for $k\not= 2$).
As is mentioned in \cite{H3}, this is an answer to a problem raised by 
Asada--Nakamura in \cite{AN}
where they pointed out the possibility that
${\mathfrak j}_{g,1}^\bQ(4)=[2]$ might be outside of $\im\tau_{g,1}^\bQ(4)$.
This possibility is now settled by Hain to be true as follows.
It is easy to see that
${\mathfrak j}_{g,1}(2)=\bZ$ and it is contained in $\im\tau_{g,1}(2)$.
However this is the whole of the
intersection of $\im\tau_{g,1}$ with the ideal ${\mathfrak j}_{g,1}$.
Namely we have
${\mathfrak j}_{g,1}(k)\cap\im\tau_{g,1}(k)=\{0\}$ for all $k\geq 3$.
Note that $\im\tau_{g,*}$ contains ${\cal L}_g$ as a natural Lie
subalgebra which corresponds to the inner derivations of
${\cal L}_g$, as was already mentioned in \S 5.
\par
In our paper \cite{Mrt9}, we introduced
certain $Sp$--equivariant mappings
$$
\Tr(2k+1)\co \h1g(2k+1)\lra S^{2k+1} H\quad (k=1,2,\cdots)
$$
such that they are surjective (after tensoring with $\bQ$)
and vanish on $\im\tau_{g,1}$ as well as on $[\h1g, \h1g]$. 
We can also show that it vanishes on the ideal
${\mathfrak j}_{g,1}$. Thus the traces descend to
$Sp$--equivariant mappings
\begin{equation*}
\Tr(2k+1)\co \hg(2k+1)\lra S^{2k+1} H\quad (k=1,2,\cdots).
\end{equation*}
These mappings are called {\it traces} because they were defined
using the trace of some matrix representation of $\h1g(2k+1)$. 
\par
We can summarize the above as follows.
The ideal ${\mathfrak j}_{g,1}$ and the image $\im\tau_{g,1}$ of the
Johnson homomorphism are almost disjoint from each other and both are
contained in the kernel of the traces.
Thus the next problem is to determine how large is the remaining
part of the Lie algebra $\h1g$.
\par
After we had found the traces, 
Nakamura \cite{N1} discovered another obstruction
to the surjectivity of $\tau_{g,1}$
which has its origin in number theory.
More precisely, it is related to
theories of outer Galois representations
of the absolute Galois group over $\bQ$ 
initiated by Grothendieck, Ihara and Deligne 
and developed by many people 
(see \cite{Gt, Ih2, D, Dr, N2} and references in them). 
Roughly speaking, Nakamura proved that certain graded
Lie algebra, which arises in the above theories,
appears in the cokernel of Johnson homomorphisms.
In particular, he concluded that the dimension of the cokernel of 
the Johnson homomorphism
$$
\tau^\bQ_{g,*}\co (\Msg(2k)/\Msg(2k+1))\otimes\bQ\lra {\mathfrak 
h}^\bQ_{g,*}(2k)
$$
is greater than or equal to a number $r_k$
which is the free $\bZ_\ell$--rank of a
certain Galois group $\Gal(\bQ(k+1)/\bQ(k))$
associated to the outer representation of
$\Gal(\overline\bQ/\bQ)$ on some nilpotent quotient
of the pro-$\ell$ fundamental
group of ${\mathbf P}^1\setminus\{0,1,\infty\}$ (see \cite{N1} for details).
Thus it may depend on the prime $\ell$. 
However it is conjectured that $r_k$ is the dimension
of the degree $k$ part of the
free graded Lie algebra over $\bQ$ which is
generated by one free generator
$\sigma_k$ of degree $k$ for each odd $k>1$ (Deligne's motivic conjecture).
As Nakamura mentions in his paper cited above, the appearance
of these Galois obstructions, in its primitive form,
was predicted by Takayuki Oda
and some closely related work was done by M. Matsumoto.
Certain estimates of these numbers have been obtained by Ihara \cite{Ih1}, 
M. Matsumoto \cite{Ma}, Tsunogai \cite{Ts} and others. 
\par
Next we consider
the images of the Johnson homomorphisms 
$$
\tau_g(k)\co\Mg(k)/\Mg(k+1)\ra \hg(k)
$$
for low degrees. The following results have been obtained. 
$$
\begin{aligned}
\im\tau_g^\bQ(1)&=[1^3]\\
\im\tau_g^\bQ(2)&=[2^2] \\
\im\tau_g^\bQ(3)&=[31^2]
\end{aligned}
\begin{aligned}
&\quad\text{(\cite{Jo2})}\\
&\quad\text{(\cite{Mrt5, H3})}\\
&\quad\text{(\cite{AN, H3})}
\end{aligned}
$$ 
Nakamura made a list of irreducible decompositions of the
$Sp$--modules ${\mathfrak h}^\bQ_{g,1}(k)$ for small $k$ by a computer 
calculation.
Utilizing it, we made 
rather long explicit computations of the Lie
bracket $\im\tau^\bQ_{g,1}(1)\otimes
\im\tau^\bQ_{g,1}(3)\ra{\mathfrak h}^\bQ_{g,1}(4)$,
in the framework of symplectic representation theory,
and obtained the following result.

\begin{proposition}
We have
$$
\im\tau_g^\bQ(4)=[42]+[31^3]+[2^3]+[31]+[2]\quad (g\geq 4)
$$
and the cokernel of the homomorphism 
$\tau_{g,1}^\bQ(4)\co (\M1g(4)/\M1g(5))\otimes\bQ\ra {\mathfrak h}^\bQ_{g,1}(4)$
is isomorphic to $[2]+[21^2]$ where $[2]={\mathfrak j}_{g,1}^\bQ(4)$.
\qed
\end{proposition}

The following list describes the structure of ${\mathfrak h}^\bQ_{g,1}(k)$
for degrees $k\leq 4$. 

$$
\vbox{
\halign{\hfil#\hfil&&\quad\hfil#\hfil\cr
$k$ & 
${\mathfrak j}_{g,1}^\bQ(k)$ & 
${\cal L}_{g}^\bQ(k)$ &
$\im\tau_g^\bQ(k)$ & 
$\Cok\tau_g^\bQ(k)\setminus\Tr$ & 
$\Tr$ \cr
\noalign{\smallskip}
\noalign{\hrule}
\noalign{\medskip}
1 & \hfil {} & [1] & $[1^3]$ & {} & {}  \cr
2 & \hfil [0] & $[1^2]$ & $[2^2]$ & {} & {}  \cr
3 & \hfil {} & [21] & $[31^2]$ & {} & [3]  \cr
4 & \hfil [2] & $[31][21^2][2]$ 
& $[42][31^3][2^3][31][2]$ & $[21^2]$ & {}  \cr
}}
$$

Here the term $\Cok\tau_g^\bQ(k)\setminus\Tr$ means that
we exclude the trace component from the 
cokernel of the homomorphism $\tau_g^\bQ(k)$.
Thus the summand $[21^2]$ in ${\mathfrak h}^\bQ_{g,1}(4)$ is not contained in
$\im\tau_g^\bQ(4)$ and it should be considered as
a new type of obstruction for the surjectivity of $\tau_{g,1}^\bQ$
other than the ideal ${\mathfrak j}_{g,1}$, the traces and the
Golois obstructions. In the next degree 5, we found yet another new
obstruction. More precisely,
at least
one copy of $[1^3]\subset {\mathfrak h}^\bQ_{g,1}(5)$ cannot be hit by 
$\tau_{g,1}^\bQ(5)$
and this leads us to the non-triviality of $H_4(\Out F_4;\bQ)$
as will be shown in \S 6.5.
\par
Thus the cokernels of the Johnson homomorphisms seem to grow rapidly
according as the degree increases and it will require many more 
studies before we can figure out the precise form of
$\im\tau^\bQ_{g,1}$ inside ${\mathfrak h}^\bQ_{g,1}$. 
We mention that Kontsevich
gave the irreducible decomposition of ${\mathfrak h}^\bQ_{g,1}$ 
as an $Sp$--module in \cite{Ko2, Ko3}
and it would be desirable to identify $\im\tau^\bQ_{g,1}$
as an explicit $Sp$--submodule of it.
\par

\subsection{Description of $Sp$--invariant part of ${\mathfrak h}^\bQ_{g,1}(6)$}

In this subsection, we describe the $Sp$--invariant part 
${\mathfrak h}^\bQ_{g,1}(6)^{Sp}$
of the degree 6 summand of the graded Lie algebra ${\mathfrak h}^\bQ_{g,1}$
to illustrate our method.
Here we only give an outline of our computation
rather than the details.
\par
We have an exact sequence
$$
0\lra \h1g(6)\lra H\otimes{\cal L}_{g,1}(7)\lra {\cal L}_{g,1}(8)\lra 0.
$$
By Proposition \ref{prop:ec}, we know that 
$\dim (H_\bQ\otimes {\cal L}_{g,1}^\bQ(7))^{Sp}=15$ and it can be shown that
$\dim {\cal L}_{g,1}^\bQ(8)^{Sp}=10$.
Hence $\dim{\mathfrak h}^\bQ_{g,1}(6)^{Sp}=5$ in the stable range.
Also it can be shown that $\dim{\mathfrak j}^\bQ_{g,1}(6)^{Sp}=2$.
Alternatively, we can use the method described  
after Proposition \ref{prop:ec} as follows.
The set of chord diagrams with $4$ chords
has $18$ elements.
If we divide the vector space 
${\cal G}_4{\cal D}^c$ spanned by it by the (4T) relation,
then the dimension reduces to $6$ (see \cite{BN}). 
For our purpose, instead of (4T) relation, we have to put the relation 
coming from condition (i) of Proposition \ref{prop:hg}. Then we find by
an explicit computation that
$$
\dim {\mathfrak h}^\bQ_{g,1}(6)^{Sp}=
\cases
5 &\quad (g\geq 3)\\
4 &\quad (g=2)\\
1 &\quad (g=1).
\endcases
$$ 
We can also give a basis $\{\al_i;i=1,\cdots, 5\}$
of $\Hom({\mathfrak h}^\bQ_{g,1}(6),\bQ)^{Sp}$ 
in terms of linear combinations of chord diagrams.
Now we consider the following 5 elements 
$\{\xi_i;i=1,\cdots,5\}$ in
${\mathfrak h}^\bQ_{g,1}(6)^{Sp}$. We choose the first two elements 
$\xi_1,\xi_2$
to be a basis of ${\mathfrak j}^\bQ_{g,1}(6)^{Sp}$. 
On the other hand, we know that
${\mathfrak h}^\bQ_{g,1}(3)=[21]+[31^2]+[3]$
where $[21]={\cal L}_g^\bQ(3)$, $[31^2]=\im\tau_g^\bQ(3)$ 
and $[3]$ is the trace component (see the table in \S 6.1).
It turns out that the unique trivial summand
$[0]$ in each of $\La^2[21], \La^2[31^2]$ and $\La^2 [3]$
survives in ${\mathfrak h}^\bQ_{g,1}(6)^{Sp}$ under the bracket operation
$\La^2{\mathfrak h}^\bQ_{g,1}(3)\ra {\mathfrak h}^\bQ_{g,1}(6)$.
We set $\xi_3,\xi_4,\xi_5$ to be the images in
${\mathfrak h}^\bQ_{g,1}(6)^{Sp}$ of these elements in the above order.
We see that $\xi_3$ is a generator of
${\cal L}^\bQ_g(6)^{Sp}\cong\bQ$.
Clearly $\xi_4$ is contained in $\im\tau_g^\bQ(6)^{Sp}$.
Also it turns out that we can construct the first two
elements $\xi_1,\xi_2$ by 
taking brackets of suitable elements
from ${\mathfrak h}^\bQ_{g,1}(1)$ and $[3]\subset{\mathfrak h}^\bQ_{g,1}(3)$.
This is one of the supporting pieces of evidence for 
Conjecture \ref{conj:abel} given in \S 6.5 below.
Now explicit computation shows that 
$\det(\al_i(\xi_j))\not=0$.
Hence $\{\xi_i\}$ forms a basis of ${\mathfrak h}^\bQ_{g,1}(6)^{Sp}$.
\par
Summing up, we have the following table for the dimensions
of $Sp$--invariant part of ${\mathfrak h}^\bQ_{g,1}(6)$
corresponding to the ideal ${\mathfrak j}^\bQ_{g,1}$, 
inner derivations ${\cal L}^\bQ_g$,
Johnson image
$\im\tau^\bQ_g$ and the remaining part 
$\Cok\tau^\bQ_g$ (for $g\geq 3$). 
\par  
$$
\vbox{
\halign{\hfil#\hfil&&\quad\hfil#\hfil\cr
${\mathfrak j}^\bQ_{g,1}(6)^{Sp}$ & ${\cal L}^\bQ_g(6)^{Sp}$ &
$\im\tau^\bQ_g(6)^{Sp}$ & 
$\Cok\tau^\bQ_g(6)^{Sp}$ &
total \cr
\noalign{\smallskip}
\noalign{\hrule}
\noalign{\medskip}
2 &  1 & 1 & 1 & 5 \cr
}}
$$
\par

We can also see from the above table
that the dimension of $\Cok\tau^\bQ_g(6)$ is exactly equal
to one. Hence it should be equal to Nakamura's 
Galois obstruction given in \cite{N1} which is in fact 
independent of $\ell$.
\par
The way of degeneration is interesting also here. More precisely, 
if we set $g=1$, then it turns out that
any of the three elements $\xi_3,\xi_4,\xi_5$ goes
to a unique non-trivial element in ${\mathfrak h}^\bQ_{1,1}(6)^{Sp}\cong\bQ$
(up to non-zero scalars).
Observe here that if $g=1$, then the Torelli group is trivial,
the ideal ${\mathfrak j}_{1,1}$ coincides with the whole of
${\mathfrak h}_{1,1}$ and the traces are all trivial.
However the forgetful mapping $\h1g\ra {\mathfrak h}_{1,1}$ is
{\it not} a Lie algebra homomorphism so that there
is no contradiction here. It may be said that,
in genus one, topology disappears and only arithmetic remains.
\par
The above is only a special case of degree $6$ summand. 
However we can see already
here a glimpse of some general phenomena. Namely 
the elements $\xi_3,\xi_4,\xi_5$ can be defined for {\it all}
degrees $4k+2$ and turn out to be non-trivial in 
${\mathfrak h}^\bQ_{g,1}(4k+2)^{Sp}$.
For example $\xi_4$ is defined to be the image of
the $Sp$--invariant part of $\La^2 [2k+1,1^2]$ in $\im\tau^\bQ_{g,1}(4k+2)$
where $[2k+1,1^2]\subset \im\tau^\bQ_{g,1}(2k+1)$ is the summand
given in \cite{AN}.
The elements $\xi_5$ for all $k>1$ were constructed in a joint work 
with Nakamura in which we are trying to understand
the Galois obstructions topologically. 

\subsection{Kernel of ${\mathfrak t}_g\ra\im\tau^\bQ_g$ and 
invariants of $3$--manifolds}

As is well known, there are close connections between the mapping
class group and $3$ dimensional manifolds.
More precisely, the Heegaard decomposition of $3$--manifolds
gives rise to a direct correspondence between 
the two objects and the construction
of surface bundles over the circle with a given monodromy
from the mapping class group yields another such relation.
Also the genus $1$ mapping class group plays a crucial
role in the theory of Dehn surgery along framed links in  
$3$--manifolds, particularly in $S^3$ by virtue of
Kirby's fundamental result \cite{Ki}. 
Hence any invariant of $3$--manifolds naturally has
certain effect on the structure of the mapping class group.
\par
Until the discovery of the Casson invariant for homology
$3$--spheres in $1985$,
the Rohlin invariant was almost the unique invariant
for $3$--manifolds. Motivated mainly by Witten's influential work
in \cite{Wi1}, the theory of topological
invariants of $3$--manifolds has been
continually and rapidly developing. See 
\cite{RT, Koh, KiM, Muh, Oh2, LMO} as well as
their references.
It is beyond the scope of this article to review these developments.
In the following, we would like to focus on those results which
have direct relation with the structure of the mapping class group. 
\par
If we are given a 3--manifold $M$ together
with an embedded oriented surface $\Sg\subset M$, then we can associate a
new manifold $M_\varph$ to each element $\varph\in\Mg$ by
cutting $M$ along $\Sg$ and pasting it back together by $\varph$. 
In the cases where $M$ is an integral homology sphere and $\varph$
belongs to the Torelli group $\Ig$, the resultant manifold $M_\varph$
is again an integral homology sphere. Hence, given any topological invariant
$\al$ of homology 3--spheres, with values in a module $A$, 
we can define a mapping $\al\co \Ig\ra A$ by setting
$\al(\varph)=\al(M_\varph)$.
Birman and Craggs \cite{BC} first studied such mappings for the case of
the Rohlin invariant $\mu$. Johnson \cite{Jo1} extended their results to obtain 
a
complete enumeration of so-called Birman--Craggs homomorphisms.
This result played an important role in his determination of the
abelianization of the Torelli group in \cite{Jo6}.
In our papers \cite{Mrt5, Mrt6}, we studied the case of the Casson 
invariant
$\la$ and in particular we obtained an interpretation of $\la$ in
terms of the secondary characteristic classes of surface bundles
(see \cite{Mrt13} and \S 6.6 below).
This work has been generalized in two different ways. 
One is due to Lescop \cite{Les} where she obtained, among other things,
a closed formula 
which expresses how the Casson--Walker--Lescop invariant 
behaves under the cut and paste operation on $3$--manifolds.
The other is given by a series of works of Garoufalidis and Levine
\cite{GL1, GL2} where they generalized our results cited above
extensively. More precisely, we considered the effect of
Casson invariant on the structure of the Torelli group
$\Ig$ while they considered all of the finite type invariants
of homology $3$--spheres introduced by Ohtsuki \cite{Oh1, Oh2}.
In particular, they proved that Ohtsuki's filtration of the
space of homology $3$--spheres can be described in terms of
the lower central series of the Torelli group. As a corollary to
this statement, they proved that any primitive
type $3k$ invariant $\al$
gives rise to a non-trivial homomorphism
$$
\al\co {\mathfrak t}_g(2k)\lra \bQ.
$$
\par
Recall from \S 5 that we have a series of homomorphisms
$$
{\mathfrak t}_g(k)\lra {\mathfrak u}_g(k)\lra \im\tau^\bQ_g(k)
$$
and we know by Hain \cite{H1} that ${\mathfrak t}_g(k)\cong{\mathfrak u}_g(k)$
for all $k\not= 2$ and that ${\mathfrak t}_g(2)={\mathfrak u}_g(2)+[0]$.
Thus the main problem concerning the above series is the following.

\begin{problem}
Determine whether the homomorphism 
${\mathfrak t}_g(k)\ra\im\tau^\bQ_g(k)$ is injective (and hence
an isomorphism)  for $k\not=2$ or not.
\label{prob:kernel}
\end{problem}

The results of Garoufalidis and Levine mentioned above show that
the above problem is crucial also from the point of view
of the theory of invariants of $3$--manifolds.
We mention that, extending earlier results of Johnson \cite{Jo4},
Kitano \cite{Kit} proved that the $k$-th Johnson
homomorphism 
$\tau_{g,1}(k)\co\M1g(k)\ra\h1g(k)$
exactly mesures the higher Massey products of
mapping tori which are associated to elements of the subgroup $\M1g(k)$.
Hence if Problem \ref{prob:kernel} will be affirmatively solved, it would imply 
that
the finite type invariants can be described in terms of the original
Casson invariant together with the Massey products.
This might sound unlikely from the point of view of $3$--manifolds
invariants.
Also the author learned from J Murakami that the restriction,
to the Torelli group,
of the projective representation of $\Mg$ associated to
the LMO invariant given in \cite{LMO, Muj} gives rise to a unipotent 
representation (after suitable truncations).
Hence it should be described by the Torelli Lie algebra ${\mathfrak t}_g$.
This also increases the importance of Problem \ref{prob:kernel}.
\par
Although Hain \cite{H3} obtained a presentation of the Torelli
Lie algebra ${\mathfrak t}_g$ (see \S 5), it is by no means easy to determine
${\mathfrak t}_g(k)\ (k=1,2,\cdots)$ explicitly by using it.
One way to compute them is to apply 
Sullivan's theory \cite{Su2} of $1$--minimal models to the Torelli group $\Ig$
which can be described as follows. 
First of all we know by Johnson \cite{Jo2} that ${\mathfrak 
t}_g(1)=U_\bQ=[1^3]$.
It was a consequence of the results of \cite{Mrt5, Mrt6} that 
${\mathfrak t}_g(2)$
contains at least $[2^2]+[0]$ where the trivial summand $[0]$
reflects the influence of the Casson invariant on the structure
of the Torelli group.
We have an isomorphism
$$
{\mathfrak t}_g(2)^*\cong\Ker(\La^2{\mathfrak t}_g(1)^*\ra H^2(\Ig)).
$$
As was already mentioned in \S 5 (Proposition 5.2), 
Hain \cite{H3} determined the right hand side
to be precisely equal to $[2^2]+[0]$ and he concluded that 
$$
{\mathfrak t}_g(2)=[2^2]+[0]\quad (g\geq 3).
$$
\par
The next case, namely the case of degree $3$ is 
given by the following result.

\begin{proposition}
We have isomorphisms
$$
{\mathfrak t}_g(3)\cong {\mathfrak u}_g(3)\cong \im\tau_g(3)=[31^2]\quad (g\geq 
6).
$$
\end{proposition}

\begin{proof}[Sketch of Proof]
Sullivan's theory \cite{Su2} implies that 
${\mathfrak t}_g(3)$ can be identified with
the kernel of the following homomorphism
$$
\Ker\left\{
\left(
\gathered
{\mathfrak t}_g(1)^*\otimes{\mathfrak t}_g(2)^*\\
\La^2{\mathfrak t}_g(2)^*
\endgathered
\right)
\overset{d}{\lra}
\left(
\gathered
\La^3 {\mathfrak t}_g(1)^*\\
\La^2{\mathfrak t}_g(1)^*\otimes{\mathfrak t}_g(2)^*
\endgathered
\right)
\right\}
\lra H^2(\Ig)
$$
where $d\equiv 0$ on ${\mathfrak t}_g(1)^*$ and 
$d\co {\mathfrak t}_g(2)^*\ra \La^2{\mathfrak t}_g(1)^*$
is the natural injection.
The differential
$d\co\La^2{\mathfrak t}_g(2)^*\ra \La^2{\mathfrak t}_g(1)^*\otimes{\mathfrak 
t}_g(2)^*$
is given by 
$d(\al\land\beta)=d\al\otimes \beta-d\beta\otimes 
\al\ (\al,\beta\in t_g(2)^*)$.
It is easy to deduce from this fact that it is injective.
Therefore ${\mathfrak t}_g(3)^*$ is isomorphic to the kernel of
the following mapping
$$
\Ker({\mathfrak t}_g(1)^*\otimes{\mathfrak t}_g(2)^*\ra \La^3 {\mathfrak 
t}_g(1)^*)
\lra H^2(\Ig)
$$
which is given by the Massey triple product.
Hain proved in \cite{H3} that all of the higher Massey products of $\Ig$
vanishes for $g\geq 6$. 
Hence passsing to the dual, we see that ${\mathfrak t}_g(3)$ is
determined by the following exact sequence
$$
\La^3 {\mathfrak t}_g(1)\lra {\mathfrak t}_g(1)\otimes{\mathfrak t}_g(2)
\overset{[\ ,\ ]}{\lra} {\mathfrak t}_g(3)\lra 0.
$$
Here the first mapping is given by
$$
\La^3 {\mathfrak t}_g(1)\ni u\land v\land w\mapsto 
u\otimes [v,w]+v\otimes [w,u]+w\otimes [u,v]
\quad (u,v,w\in {\mathfrak t}_g(1))
$$
and the image of which in ${\mathfrak t}_g(3)$ under the bracket
operation 
is trivial because of the Jacobi identity.
The irreducible decomposition of 
${\mathfrak t}_g(1)\otimes{\mathfrak t}_g(2)=[1^3]\otimes ([2^2]+[0])$
is given by
$$
{\mathfrak t}_g(1)\otimes{\mathfrak t}_g(2)=
([3^21]+[321^2]+[2^21^3]+[31^2]+[32]+[21^3]+[2^21]
+[21]+[1^3])+[1^3].
$$
Explicit computations, corresponding to each summand
in the above decomposition,
show that the mapping 
$\La^3{\mathfrak t}_g(1)\ra{\mathfrak t}_g(1)\otimes{\mathfrak t}_g(2)$
hits any summand except $[31^2]$.
Since we already know that $\im\tau^\bQ_g(3)=[31^2]$, the result follows.
\end{proof}

If we inspect the above proof carefully,
we find that it was not necessary to use 
Hain's vanishing of Massey triple products of
$\Ig$ for $g\geq 6$. It turns out that 
the computation itself contains a proof of
the vanishing of them.
However for higher degrees $k=4,5,\cdots$, 
Hain's result considerably simplifies the computations.
By using this method, 
we can continue the computation for higher degrees. 
For example the degree $4$
summand ${\mathfrak t}_g(4)$ can be determined by the following exact sequence
$$
\La^2{\mathfrak t}_g(1)\otimes {\mathfrak t}_g(2)\lra 
\bigl({\mathfrak t}_g(1)\otimes{\mathfrak t}_g(3)\bigr)\oplus\La^2{\mathfrak 
t}_g(2)
\overset{[\ ,\ ]}{\lra} {\mathfrak t}_g(4)
\lra 0.
$$
Here the first mapping is given by
$$
(u\land v)\otimes w\mapsto u\otimes [v,w]-v\otimes [u,w]-[u,v]\land w
\in \bigl({\mathfrak t}_g(1)\otimes{\mathfrak t}_g(3)\bigr)\oplus\La^2{\mathfrak 
t}_g(2)
$$
where $(u,v\in {\mathfrak t}_g(1), w\in{\mathfrak t}_g(2))$ and the above 
element
vanishes in ${\mathfrak t}_g(4)$ again by the Jacobi identity.
Although our computation is not finished yet, we see no signs of
non-trivial kernel for 
${\mathfrak t}_g(k)\ra\im\tau^\bQ_g(k)\ (k=4,5,6)$ 
so far.
\par
We can also ask how other invariants of homology $3$--spheres,
for example various Betti numbers of Floer homology \cite{Fl} or infinitely many 
homology cobordism invariants the existence of which
is guaranteed by a remarkable result
of Furuta \cite{Fu}, will influence the structure of the Torelli group.
Also it seems to be a challenging problem to seek those invariants
of homology $3$--spheres which reflect semi-simple informations,
rather than the nilpotent ones,
of the Torelli group.

\subsection{Cocycles for the Mumford--Morita--Miller 
classes}

In our paper \cite{Mrt10}, we constructed certain representations
$\rho_1$ of the mapping class groups $\Msg, \Mg$ and obtained
the following commutative diagram
\begin{equation*}
\begin{CD}
\Msg @>{\rho_1}>> \frac{1}{2}\La^3 H\rtimes Sp(2g, \bZ)\\
@VVV @VVV\\
\Mg @>>{\rho_1}> \frac{1}{2}\La^3 H/H\rtimes Sp(2g, \bZ).
\end{CD}
\end{equation*}

By making use of a standard fact concerning the cohomology of
a group which is a semi-direct product, we deduced the existence
of the following diagram
\begin{equation}
\begin{CD}
\Hom(\La^*(\La^3 H_\bQ), \bQ)^{Sp} @>{\rho_1^*}>> H^*(\Msg;\bQ)\\
@AAA @AAA \\
\Hom(\La^* U_\bQ, \bQ)^{Sp} @>>{\rho_1^*}>  H^*(\Mg;\bQ).
\end{CD}
\label{7}
\end{equation}

If we combine this with Proposition \ref{prop:alpha} and
Proposition \ref{prop:beta} (see \S 4.2), then we obtain
the following commutative diagram (which is defined
at the cocycle level)
\begin{equation}
\begin{CD}
\bQ[\al_\vGa;\vGa\in{\cal G}] @>{\rho_1^*} >> H^*(\Msg;\bQ)\\
@AAA @AAA \\
\bQ[\beta_\vGa;\vGa\in{\cal G}^0] @>>{\rho_1^*}> H^*(\Mg;\bQ).
\end{CD}
\label{8}
\end{equation}

It was proved in \cite{Mrt12} 
that the images of $\rho_1^*$ in \thetag{8} contain
all of the characteristic classes $e, e_i\ (i=1,2,\cdots)$. 
However the problem of determining 
whether the images contain new classes or not remained open.
This problem was soon solved by the introduction of 
generalized Mumford--Morita--Miller classes due to Kawazumi 
in \cite{Kaw5}.
There is a remarkable work of Looijenga \cite{Lo4}
in which he determined the stable cohomology of
$\Mg$ with coefficients in any finite dimensional irreducible
representation of $Sp(2g,\bQ)$.
In \cite{Kaw4}, Kawazumi used the above generalized classes
to give a different basis for some of these  
twisted cohomology groups,
thereby adding a topological flavor to Looijenga's result.
We also mention that Ivanov \cite{Iv} obtained a stability theorem
for twisted cohomology groups of mapping class groups
(for a surface with at least one boundary component).
\par
Let us define ${\cal R}^*(\Mg)$ to be the
subalgebra of $H^*(\Mg;\bQ)$ generated by the classes $e_i$. 
Similarly we define ${\cal R}^*(\Msg)$ to be the subalgebra
of $H^*(\Msg;\bQ)$ generated by the classes $e, e_i$. We may call
them the {\it tautological algebra} of the mapping class groups.
These are just the reduction to the rational cohomology
of the original tautological algebras
${\cal R}^*({\mathbf M}_g), {\cal R}^*({\mathbf C}_g)$ of the moduli spaces
which are defined to be subalgebras of the individual Chow algebras
generated by the classes $\kappa_i$ and $c_1(\om)$ 
(see \cite{Fa4, HL, HM}). 

\begin{theorem} [Kawazumi--Morita \cite{KM1, KM2}]
For any $g$,
the images of $\rho_1^*$ in \thetag{8} coincide with the tautological
algebras ${\cal R}^*(\Msg)$ and ${\cal R}^*(\Mg)$ of the mapping class groups.
\qed
\label{thm:KM}
\end{theorem}
\par
Thus the homomorphisms $\rho_1^*$ in \eqref{7} have rather big kernel.
For the lower homomorphism $\rho_1^*$,
Garoufalidis and Nakamura \cite{GN} have given an interpretation of this
fact in the framework of symplectic representation theory
by showing an isomorphism
$$
\bigl(\La^* U_\bQ^*/([2^2])\bigr)^{Sp}\cong \bQ[e_1,e_2,\cdots]
$$
which holds in the stable range, where $([2^2])$ denotes the ideal of
$\La^*U_\bQ^*$ generated by $[2^2]\subset \La^2U_\bQ^*$.
Their result can be generalized to the case of the upper $\rho_1^*$
so that we have an isomorphism
$$
\bigr(\La^*(\La^3 H_\bQ^*)/([2^2]+[1^2])\bigr)^{Sp}\cong \bQ[e,e_1,e_2,\cdots]
$$
which holds also in the stable range. Here $[1^2]$ is a
certain {\it diagonal} summand in
$\La^2(\La^3 H_\bQ^*)$, which has three copies of $[1^2]$,
described explicitly in \cite{Mrt12}. 
However these results are, at present, valid only in the stable
range while Theorem \ref{thm:KM} is true in all degrees.
If we pass to the dual context, namely if we consider the homomorphism
$$
(\rho_1)_*\co H_*(\Mg;\bQ)\lra 
H_*(\frac{1}{2}\La^3 H/H\rtimes Sp(2g, \bZ);\bQ)
$$
induced by $\rho_1$ on homology, then we find that
those cycles in $(\La^* U_\bQ)^{Sp}$ which come from the moduli space
must have fairly restricted types.
In particular, Faber's result below (Theorem \ref{thm:Fa}) implies that
the subspace of $(\La^{2g-4}U_\bQ)^{Sp}$  
generated by the
{\it moduli cycles} is one dimensional and
a similar statement is valid for $(\La^{2g-2}(\La^3 H_\bQ))^{Sp}$.
Hence we can define {\it fundamental cycles}
$$
\mu_{g,*}\in(\La^{2g-2}(\La^3 H_\bQ))^{Sp},\quad
\mu_g\in (\La^{2g-4}U_\bQ)^{Sp}
$$
which are well defined up to scalars and which should be
expressed in terms of certain linear combinations of
$Sp$--invariant tensors $a_\vGa, b_\vGa$ described in \S 4.
\par
Starting from basic works on the Chow algebras of the
moduli spaces ${\mathbf M}_3, {\mathbf M}_4$ in \cite{Fa1, Fa2}, Faber
made numerous explicit computations concerning the tautological
algebra of the moduli spaces.
Based on them, he  
proposed the following beautiful conjecture about the structure of the
tautological algebra ${\cal R}^*({\mathbf M}_g)$ (see \cite{Fa4} for details, 
in particular for the precise description of part (3) below). 

\begin{conjecture}{\rm (Faber \cite{Fa4})}
\begin{enumerate}
\item The tautological algebra ${\cal R}^*({\mathbf M}_g)$ of the
moduli space ${\mathbf M}_g$ ``behaves like" the 
cohomology algebra of a nonsingular projective variety of dimension
$g-2$. More precisely, it vanishes in degrees $> g-2$, is one dimensional
in degree $g-2$ and the natural pairing ${\cal R}^{i}({\mathbf M}_g)\times
{\cal R}^{g-2-i}({\mathbf M}_g)\ra {\cal R}^{g-2}({\mathbf M}_g)$ is perfect.
It also satisfies the Hard Lefschetz and the Hodge Positivity properties
with respect to the class $\kappa_1$.
\item The $[\frac{g}{3}]$ classes $\kappa_1,\cdots,\kappa_{[\frac{g}{3}]}$
generate the algebra with no relations in degrees $\leq [\frac{g}{3}]$.
\item There exist explicit formulas for the proportionalities in
degree $g-2$.
\end{enumerate}
\label{conj:Faber}
\end{conjecture}

Several supporting pieces of evidence for this conjecture have been obtained. 
\par
\begin{theorem}[Looijenga \cite{Lo3}]
The tautological algebra ${\cal R}^*({\mathbf M}_g)$ is trivial in degrees
$> g-2$ and ${\cal R}^{g-2}({\mathbf M}_g)$ is at most one dimensional.  
Similarly 
${\cal R}^*({\mathbf C}_g)$ is trivial for $*>g-1$
and ${\cal R}^{g-1}({\mathbf C}_g)$ is at most one dimensional.
\qed
\end{theorem}
\par
We mention that Jekel \cite{Je} proved the vanishing 
$e^{g}=0\in H^{2g}(\Msg;\bQ)$ by making use of
a certain representation $\Msg\ra \ho_+ S^1$
(cf \cite{Mrt3}). This
gives a purely topological proof of a part of Looijenga's
result above. 
Using the results of Mumford in \cite{Mu} as well as
those of Witten \cite{Wi2} and Kontsevich \cite{Ko1}
combined with Looijenga's theorem above, Faber proved the
following.

\begin{theorem}[Faber \cite{Fa3}]
$\kappa_{g-2}$ is non-zero on ${\mathbf M}_g$ so that 
${\cal R}^{g-2}({\mathbf M}_g)$
is one dimensional. 
\label{thm:Fa}
\qed
\end{theorem}

It follows immediately that ${\cal R}^{g-1}({\mathbf C}_g)$ is 
also one dimensional.
\par
The above results are obtained mainly in the framework of 
algebraic geometry.
Here we would like to describe a topological approach to Faber's
conjecture which has fairly different feature. Naturally we can 
obtain relations only in the rational cohomology algebra rather
than the Chow algebra. However we hope that this approach would
have its own meaning.
\par
Our method is very simple. Namely, 
associated to any unstable relation in $H^*(\La^3 H_\bQ)^{Sp}$ or
in $H^*(U_\bQ)^{Sp}$,
we can obtain a polynomial
relation in the tautological algebra 
${\cal R}^*(\Msg)$ or ${\cal R}^*(\Mg)$ by applying Theorem \ref{thm:KM}.
For example, we can apply Proposition \ref{prop:relation} 
to obtain a series
of non-trivial
relations as follows. We know by this proposition that
there is a unique relation

\begin{equation}
\sum_{C\in{\cal D}^\ell(6k)}\ a_C=0
\label{9}
\end{equation}

in $(H_\bQ^{6k})^{Sp}$ for $g=3k-1$. Fortunately,
this relation survives
in $(\La^{2k} U_\bQ)^{Sp}$ under the natural projection
$H_\bQ^{6k}\ra \La^{2k} U_\bQ$. Hence passing to the dual,
the relation $\sum_C \al_C=0$ gives rise to a polynomial relation
in ${\cal R}^{2k}({\cal M}_{3k-1})$ which expresses the class
$e_{k}$ as a polynomial in lower $e_i$'s.
It turns out that this relation is exactly the same as Faber's
relation mentioned in \cite{Fa4} up to a factor of some powers of $2$.
The associated generating function appears, in our context,
as a result of enumeration of certain trivalent graphs.
One of the merits of our method is that once we obtain a relation
in ${\cal R}^*(\Mg)$ for some $g$, we can obtain associated relations
for all genera $< g$. This is because of the following reason.
Although the mapping
$H_1(\Sigma_{g};\bQ)\ra H_1(\Sigma_{g-1};\bQ)$, which is induced
by collapsing the last handle, is not very natural from the point of view
of algebraic geometry, it does induce a natural mapping
$$
(H_g^{\otimes 2k})^{Sp(2g,\bQ)} \lra 
(H_{g-1}^{\otimes 2k})^{Sp(2g-2,\bQ)}
$$
where $H_g$ and $H_{g-1}$ stand for 
$H_1(\Sigma_{g};\bQ)$ and $H_1(\Sigma_{g-1};\bQ)$
respectively. For example, the relation \eqref{9}
which is the unique relation for $g=3k-1$ continues
to hold for all $g<3k-1$.
In this way, using the unique relation
\eqref{9} above, we can prove the following theorem. 

\begin{theorem}
The tautological algebra ${\cal R}^*(\Mg)$ of $\Mg$ is generated
by the first $[\frac{g}{3}]$ Mumford--Morita--Miller classes
$$
e_1,e_2,\cdots,e_{[\frac{g}{3}]}.
$$
Moreover there are explicit formulas which
express any class $e_j\ (j>[\frac{g}{3}])$ 
as a polynomial in the above classes.
\qed
\label{thm:Morita}
\end{theorem}

This theorem gives an affirmative solution 
(at the level of rational cohomology) to a part of Faber's
conjecture (Conjecture \ref{conj:Faber}).
Details will be given in a forthcoming paper \cite{Mrt14}.
We expect that we can obtain further relations in ${\cal R}^*(\Mg)$
as well as in ${\cal R}^*(\Msg)$
by this method. We also expect that the explicit form of
the fundamental cycle $\mu_g$ mentioned above 
should be closely related to part (3) of Conjecture \ref{conj:Faber}.
\par
Concerning the (co)homology of the moduli space ${\mathbf M}_g$,
there is another well known conjecture due to Witten and Kontsevich
(see \cite{Ko1}) which says that certain natural cycles of ${\mathbf M}_g$
constructed by them should be expressed in terms of the
Mumford--Morita--Miller classes. The first case was affirmatively
solved by Penner in \cite{Pe} and Arbarello and Cornalba made
considerable progress on this problem in \cite{AC}.
We expect that there would exist certain connection between
their method in the framework of algebraic geometry with 
our approach given above which uses symplectic representation theory.
Penner \cite{Pe2} described a related conjecture in the context
of his new model of a universal Teichm\"uller space.
\par
Also there is a problem concerning the {\it unstable} cohomology
of $\Mg$ or ${\mathbf M}_g$. Harer and Zagier pointed out in \cite{HZ} that
their determination of the orbifold Euler characteristic of ${\mathbf M}_g$
implies that there must exist many unstable cohomology classes.
However it seems that the unstable class constructed by
Looijenga \cite{Lo1} for genus $3$ moduli space
is the only known explicit example.

\begin{problem}
Construct unstable cohomology classes of $\Mg$ explicitly.
\end{problem}

\subsection{Cohomology of the graded Lie algebra ${\mathfrak h}^\bQ_\infty$ and
$\Out F_n$ }

Here we consider the (co)homology of the graded Lie algebra
$$
{\mathfrak h}_{g,1}=\bigoplus_{k=1}^\infty {\mathfrak h}_{g,1}(k).
$$
First we would like to know the abelianization of $\h1g$.
The first Johnson homomorphism together with the traces
induce a graded Lie algebra homomorphism
$$
(\tau_{g,1}(1),\oplus_k \Tr(2k+1))\co \h1g
\lra \La^3 H\oplus \bigl(\bigoplus_{k\geq 1} S^{2k+1} H\bigr)
$$
which is surjective after tensoring with $\bQ$
where the target is considered as an {\it abelian} Lie algebra.
As in \S 5, let us consider the direct limit
$$
{\mathfrak h}_\infty=\lim_{g\to\infty}\h1g,\quad 
{\mathfrak h}_\infty^\bQ=\lim_{g\to\infty}{\mathfrak h}^\bQ_{g,1}.
$$
It is easy to see that the Johnson homomorphism $\tau_{g,1}(1)$ and
the traces $\Tr(2k+1)$ are all compatible with respect to the 
inclusions $\h1g\subset {\mathfrak h}_{g+1,1}$. Namely the following
diagram is commutative
$$
\begin{CD}
\h1g @>>> \La^3 H_g\oplus \bigl(\bigoplus_{k\geq 1} S^{2k+1} H_g\bigr) \\
@VVV@VVV \\
{\mathfrak h}_{g+1,1} @>>> \La^3 H_{g+1}\oplus 
\bigl(\bigoplus_{k\geq 1} S^{2k+1} H_{g+1}\bigr) \\
\end{CD}
$$
where $H_g$ denotes the module $H$ corresponding to the genus $g$.
In view of the explicit computations in low degrees
we have done so far, it seems to be
reasonable to make the following conjecture.

\begin{conjecture}
The abelianization of the Lie algebra 
${\mathfrak h}^\bQ_{g,1}$ is given by $\tau_{g,1}(1)$ and the traces,
so that we have an isomorphism
$$
H_1({\mathfrak h}^\bQ_{g,1})=
{\mathfrak h}^\bQ_{g,1}/[{\mathfrak h}^\bQ_{g,1},{\mathfrak h}^\bQ_{g,1}]\cong 
\La^3 H_\bQ\oplus \bigl(\bigoplus_{k\geq 1} S^{2k+1} H_\bQ\bigr).
$$
A similar statement is true for ${\mathfrak h}_\infty^\bQ$.
\label{conj:abel}
\end{conjecture}

If this conjecture were true, then any element in ${\mathfrak h}^\bQ_{g,1}$,
including the Galois obstructions,
can be described by taking brackets of suitable elements of
$\La^3 H_\bQ$ and $S^{2k+1} H_\bQ$.
Regardless of whether the above conjecture is true or not, we have
a homomorphism
$$
H_c^*\bigl(\La^3 H_\bQ\oplus \bigl(\bigoplus_{k\geq 1} 
S^{2k+1} H_\bQ\bigr)\bigr)
\lra H_c^*({\mathfrak h}^\bQ_{g,1})
$$
where $H_c^*$ denotes the continuous cohomology,
in the usual sense (cf \cite{Ko2}),
of graded Lie algebras 
which are infinite dimensional but each degree $k$ summand 
is finite dimensional for all $k$ .
The $Sp$--invariant part of the above homomorphism can be written as 
\begin{equation}
\begin{split}
&H_c^*\bigl(\La^3 H_\bQ\oplus 
\bigl(\bigoplus_{k\geq 1} S^{2k+1} H_\bQ\bigr)\bigr)^{Sp}\\
\cong& \bigl(\La^*(\La^3 H^*_\bQ)\otimes
\La^* (S^3 H^*_\bQ)\otimes \La^* (S^5 H^*_\bQ)
\otimes\cdots\bigr)^{Sp}
\lra H_{c}^*({\mathfrak h}_{\infty}^{\bQ})^{Sp}.
\end{split}
\label{10}
\end{equation}
where $H_c^*({\mathfrak h}_\infty^\bQ)^{Sp}$ denotes the inverse limit
of the $Sp$--invariant part of the continuous
cohomology of the graded Lie algebras ${\mathfrak h}^\bQ_{g,1}$ 
which stabilizes.
The Lie algebra ${\mathfrak h}^\bQ_{g,1}$ contains $\im\tau^\bQ_{g,1}$ as a 
natural
Lie subalgebra and by restriction and passing to the limit, 
we obtain a series of
homomorphisms
\begin{equation}
\begin{split}
H^*_c({\mathfrak h}_\infty^\bQ)^{Sp}\lra H^*_c(\Gr {\mathfrak n}_\infty)^{Sp}
&\cong
H^*_c({\mathfrak n}_\infty)^{Sp}
\cong \lim_{g\to\infty} H^*({\cal N}_{g,1})^{Sp}\\
&\lra
\lim_{g\to\infty} H^*(\M1g).
\end{split}
\label{11}
\end{equation}

Here ${\mathfrak n}_\infty=\lim_{g\to\infty} {\mathfrak n}_{g,1}\ ({\mathfrak 
n}_{g,1}(k)
\cong \im \tau^\bQ_{g,1}(k)\subset{\mathfrak h}^\bQ_{g,1}(k))$
and the second isomorphism is due to Hain \cite{H3} where he 
showed that ${\mathfrak n}_{g,1}$ has a mixed Hodge structure (depending
on a fixed complex structure on the reference surface)
so that there is a canonical isomorphism of ${\mathfrak n}_{g,1}$ with
$\Gr{\mathfrak n}_{g,1}$ after tensoring with $\bC$.
For the last homomorphism, see \cite{Mrt12}.
Since $\Tr(2k+1)$ is trivial on $\im \tau_{g,1}(2k+1)$ 
for any $k$ (see \cite{Mrt9}),
the composition of \eqref{10} with \eqref{11} is trivial
on any $Sp$--invariant which contains the trace
component $S^{2k+1} H_\bQ^*$.
Hence, the {\it moduli part} of the homomorphism \eqref{10}
is given by 
$$
\La^*(\La^3 H_\bQ^*)^{Sp}\lra H^*_c({\mathfrak h}_\infty^\bQ)^{Sp}.
$$
As was explained in \cite{HL, KM1} (see also \S 2), 
a combination of our result
\cite{KM1} with that of \cite{H3} implies that the image of the above
homomorphism can be identified
with the polynomial algebra
$$
\bQ[e_1,e_2,\cdots]\subset \lim_{g\to\infty} H^*(\M1g;\bQ).
$$
Then it is a natural question to ask the geometric meaning of
the remaining part of the homomorphism \eqref{10}.
This can be answered by invoking important work of 
Kontsevich \cite{Ko2, Ko3}, 
which we now
briefly review. 
\par
As was mentioned in \S 5, our ${\mathfrak h}_\infty^\bQ$
is the degree positive part of his Lie algebra $\ell_\infty$
described in the above cited papers.
The homology group $H_*(\ell_\infty)$ has a natural structure of 
a commutative and cocommutative Hopf algebra, where the
multiplication comes from the sum operation 
$\ell_\infty\oplus\ell_\infty\ra\ell_\infty$. Also it can be decomposed
as
$$
H_*(\ell_\infty)\cong H_*({\mathfrak{sp}}(2\infty;\bQ))\otimes 
H_*({\mathfrak h}_\infty^\bQ)_{Sp}.
$$
By making use of his theory of graph cohomology
together with a result of Culler--Vogtmann \cite{CV}, he 
proved the following remarkable theorem.

\theorem[Kontsevich \cite{Ko2, Ko3}]
There exists an isomorphism
$$
PH_k(\ell_\infty)\cong PH_k({\mathfrak{sp}}(2\infty;\bQ))\oplus\bigoplus_{n\geq 
2}
H^{2n-2-k}(\Out F_n;\bQ)
$$
where $PH_k$ denotes the primitive part of the $k$--dimensional homology and
$\Out F_n$
denotes the outer automorphism group of the free group of rank $n$.
More precisely, for each even degree $2n\ (n>0)$ 
with respect to the grading of
$\ell_\infty$, there exists an isomorphism
$$
PH_k(\ell_\infty)_{2n}\cong
H^{2n-k}(\Out F_{n+1};\bQ).
$$
Passing to the dual, we also have an isomorphism
$$
PH^k_c({\mathfrak h}_\infty^\bQ)^{Sp}_{2n}\cong
H_{2n-k}(\Out F_{n+1};\bQ).\eqno{\qed}
$$
\label{thm:Kontsevich}
\rm

Let us observe here that if Conjecture \ref{conj:abel} were true, then
$H_1({\mathfrak h}_\infty^\bQ)_{Sp}=0$ so that we can conclude
$H^{2n-3}(\Out F_n;\bQ)=0$ for any $n\geq 2$ by the above
theorem of Kontsevich. We have checked that
this is the case up to $n=4$.
We mention that Culler--Vogtmann proved,
in the above paper, that the virtual cohomological dimension of
$\Out F_n$ is equal to $2n-3$. 
\par
Now we apply Theorem \ref{thm:Kontsevich} to the homomorphism \eqref{10}. 
We know that there is a copy $S^3 H_\bQ=[3]\subset {\mathfrak h}^\bQ_{g,1}(3)$
which goes to $S^3 H_\bQ$ bijectively by the trace $\Tr(3)$. 

\begin{proposition}
The homomorphism 
$$
\La^2 S^3 H_\bQ \ni \xi\land \eta\longmapsto [\xi,\eta]\in {\mathfrak 
h}^\bQ_{g,1}(6)
\quad (\xi, \eta\in S^3 H_\bQ)
$$
is injective.
\label{prop:la2}
\end{proposition}

\begin{proof}[Sketch of Proof]
It is easy to see that the irreducible decomposition of 
the module
$\La^2 S^3 H_\bQ=\La^2 [3]$ 
is given by
$$
\La^2 S^3 H_\bQ=[51]+[4]+[3^2]+[2^2]+[1^2]+[0].
$$
Then the result follows from rather
long explicit computations, in the framework of
symplectic representation theory, of the Lie bracket
$\La^2 [3]\ra {\mathfrak h}^\bQ_{g,1}(6)$.
\end{proof}
\par
Recall that similar homomorphism
$$
\La^2 (\La^3 H_\bQ)\lra {\mathfrak h}^\bQ_{g,1}(2)
$$
has a big kernel (Hain \cite{H3}, see also \cite{Mrt12} for the case of one 
boundary
component).
Thus we find that the behaviour of the trace component  
is fairly different
from that of the moduli part. In particular, the pull back of any class
in $H^2(S^3 H_\bQ)$ to the subalgebra of ${\mathfrak h}^\bQ_{g,1}$ 
generated by $S^3 H_\bQ$
is trivial. However if we consider the interaction of the trace
component with the moduli part, the situation changes.
In fact, we obtain the following result.

\begin{proposition}
The pull back of the $Sp$--invariant part $H^2(S^3 H_\bQ)^{Sp}\cong\bQ$
to $H^2_c({\mathfrak h}_\infty^\bQ)^{Sp}$ by the trace $\Tr(3)$ is a
non-trivial primitive element.
\end{proposition}

\begin{proof}[Sketch of Proof]
By Proposition \ref{prop:la2}, there is a chain 
$u=\sum_i (a_i,b_i)\in C_2({\mathfrak h}^\bQ_{g,1})$ such that 
$a_i,b_i\in S^3 H_\bQ\subset{\mathfrak h}^\bQ_{g,1}(3)$
and $\partial u$ is a non-zero element of ${\mathfrak h}^\bQ_{g,1}(6)^{Sp}$.
On the other hand, explicit computation shows that
there is a chain $v=\sum_j (a_j',b_j')\in C_2({\mathfrak h}^\bQ_{g,1})$
such that $a'_j\in {\mathfrak h}^\bQ_{g,1}(1), b'_j\in{\mathfrak 
h}^\bQ_{g,1}(5)$ and
$\partial v=\partial u$. We can arrange the above elements so that
the $2$--cycle $u-v$ is in fact an $Sp$--invariant one.
Thus it defines an element of $H_2({\mathfrak h}_\infty^\bQ)^{Sp}_6$
on which the cohomology class in question takes a non-zero
value. The primitivity follows from the property of the trace.
\end{proof}
\par
If we combine the above result with Thoerem \ref{thm:Kontsevich}, we can
conclude the non-triviality of $H_4(\Out F_4;\bQ)$.
This seems to be consistent with a recent result of Hatcher--Vogtmann \cite{HV}
in which they proved that $H_4(\Aut F_4;\bQ)\cong\bQ$
(Vogtmann informed us that she obtained
an isomorphism $H_4(\Aut F_4;\bQ)\cong H_4(\Out F_4;\bQ)$
by a computer calculation).
\par
If we use higher traces $\Tr(2k+1)$ in the above
consideration, we obtain 
infinitely many classes in $H^2_c({\mathfrak h}_\infty^\bQ)^{Sp}$
and these in turn give rise to
a series of certain homology classes in 
$$
H_{4k}(\Out F_{2k+2};\bQ)\quad (k=1,2,\cdots).
$$
It seems highly likely that all of these classes are non-trivial.
Also by combining various trace components as well as the moduli
part at the same time, we can define many primitive (co)homology classes of 
the Lie algebra $\ell_\infty$. 

\subsection{Secondary characteristic classes of surface bundles}

Chern classes and Pontrjagin classes are representatives
of characteristic classes of vector bundles and they play
fundamental roles in diverse branches of mathematics. We may call these
classes {\it primary} characteristic classes. In some cases
where these primary classes vanish, there arise various
theories of so-called {\it secondary} characteristic classes.
For example, we have the theory of characteristic classes of
foliations, characteristic classes of flat bundles,
the theory of Chern and Simons and also that of Cheeger and Simons
and so on. 

In the case of surface bundles, it is natural to call the classes
$e_i$ the primary characteristic classes. As was mentioned in
\S 3, the odd classes $e_{2i-1}$ come from the Siegel modular 
group $Sp(2g,\bZ)$ via the classical representation 
$\rho_0\co \Mg\ra Sp(2g,\bZ)$. Hence these classes (with rational
coefficients) vanish on the Torelli
group $\Ig$. It is a fundamental question
concerning the cohomology of $\Ig$ whether even classes 
$e_{2i}$ are non-trivial
on it or not (see Conjecture 3.4). 
On the other hand, it was proved in
\cite{Mrt12} that any class $e_i$, including the cases of even $i$,
vanish on the subgroup
$\Kg$ of the mapping class group $\Mg$, which
is the subgroup generated by all Dehn twists along {\it separating}
simple closed curves on $\Sg$ (see \S 2). 
The proof of this fact goes roughly as follows. 
As was recalled in \S 6.4, in our paper \cite{Mrt12} we have 
constructed certain natural cocycles for $e_i$ and by the very definition
they vanish on $\Kg$.
Thus there are {\it two} reasons with different sources that the
odd classes $e_{2i-1}$ vanish on $\Kg$.
\par
By making use of this fact, we can define the {\it secondary}
characteristic classes of surface bundles. One way to do so
can be described as follows. For each $e_{2i-1}$,
choose two cocycles $c, c'$ both of which represent $e_{2i-1}$ such that
$c$ comes from the Siegel modular group while $c'$ is a cocycle
constructed in \cite{Mrt12} using the linear
representation $\rho_1$ described there.
Since these two cocycles are cohomologous to each other,
there exists a cochain
$$
d\in C^{4i-3}(\Mg;\bQ)
$$
such that $\delta d=c-c'$. Now both of $c, c'$ are 0 on $\Kg$
so that if we restrict $d$ to $\Kg$, it is a cocycle.
Hence we can define a cohomology class
$$
d_i\in H^{4i-3}(\Kg;\bQ)
$$
to be the class of the cocycle $d\in Z^{4i-3}(\Kg;\bQ)$.
It can be shown that the cohomology class of $[d]$ does not
depend on the choices of the cocycles $c, c'$ modulo the
indeterminacy 
$$
\im \bigl(H^{4i-3}(\Mg;\bQ)\lra H^{4i-3}(\Kg;\bQ)\bigr).
$$
We remark that if Conjecture \ref{conj:basic} were true, then the above
indeterminacy vanishes, at least in the stable
range, so that $d_i$
would be uniquely defined.
\par
We know that the first one
$$
d_1\co\Kg\lra \bQ
$$
can be defined uniquely. Here we describe the precise formula
for it. 
Let $\tau\in Z^2(\Mg;\bZ)$ be Meyer's signature cocyle given
in \cite{Mey} which is in fact defined on $Sp(2g,\bZ)$. It represents
$-\frac{1}{3} e_1$. On the other hand in \cite{Mrt4, Mrt8} we gave another
cocycle for $e_1$. More precisely, in the notation of \S 6.4
the cocycle 
$$
c'=\frac{1}{2g+1}\ (-3 \al_{\vGa_1}+(2g-2) \al_{\vGa_2})
$$
represents $e_1$
(see \cite{HR} for an interpretation of this result
as well as others from the point of view of
algebraic geometry). Here $\vGa_1$ is the trivalent graph which
has $2$ vertices and $2$ loops while $\vGa_2$ is the 
trivalent graph with $2$ vertices and without loops (namely
the {\it theta} graph). The above two cocycles are cohomologous
to each other so that there exists a mapping
$$
d_1\co\Mg\lra \bQ\quad (g\geq 2)
$$
such that $\delta d_1=-3\tau-c'$. Since $\Mg$ is perfect for $g\geq 3$
and $H_1({\cal M}_2)=\bZ/10$, the above map $d_1$ is uniquely defined.
Also since the restriction of both of $\tau$ and $c'$ to the subgroup
$\Kg$ is trivial, we obtain a homomorphism
$$
d_1\co\Kg\lra \bQ.
$$
This is the definition of our secondary class $d_1\in H^1(\Kg;\bQ)$.

\begin{theorem}[Morita \cite{Mrt6}]
Let $\varph\in\Kg$ be a Dehn twist along a separating simple
closed curve on $\Sg$ such that it divides $\Sg$ into
two compact surfaces of genus $h$ and $g-h$. Then
the value of the secondary class $d_1\in H^1(\Kg;\bQ)$ on it
is given by
$$
d_1(\varph) = \frac{12}{2g+1}\ h(g-h).
$$
Moreover $d_1$ is the generator of $H^1(\Kg;\bQ)^{\Mg}\cong\bQ$ for
all $g\geq 2$.
\qed
\end{theorem}

In our paper \cite{Mrt6}, the coefficient of $h(g-h)$
in the above formula for $d_1(\varph)$ 
was not mentioned. However it is easy to deduce it from the
results obtained there. 
In \cite{Mrt13} we gave an interpretation of $d_1$ in terms of
Hirzebruch's signature defect of certain  framed 3--manifolds.
Generalizing this, we obtained another 
more geometrical definition of higher secondary
classes $d_i\in H^{4i-3}(\K1g;\bQ)$. We expect that these two
definitions of $d_i$ would coincide for all $i$.

\begin{remark}
Hain informed us the following interesting facts.
In the above theorem, the number $h(g-h)$ appeared in relation
to the Dehn twist $\varph$. About the same time, 
the same number appeared in
the works of Jorgenson \cite{Jor} and Wentworth \cite{We}.
In fact, it plays the role of the
principal factor in their asymptotic
formula of Faltings delta function around the divisor
$\Delta_h$, which is associated to $\varph$,
of the Deligne--Mumford compactification $\overline{\mathbf M}_g$.
This number plays one more important role also in a recent work
of Moriwaki \cite{Mrw} where he obtained a certain inequality 
related to the cone of positive divisors on $\overline{\mathbf M}_g$.
In fact, Hain has an interpretation that these phenomena are 
more than just a coincidence. 
Kawazumi is also trying to develop a theory related to $d_1$.
From the topological point of view,
we may say that they are the manifestation of the Casson invariant
in the geometry of the moduli space of curves. 
We expect that there should exist very rich structures here
which deserve future investigations.
\end{remark}

\begin{conjecture} 
All of the secondary classes $d_i\in H^{4i-3}(\Kg;\bQ)$ are
uniquely defined and non-trivial if $g$ is 
sufficiently large. Moreover we have an isomorphism
$$
\lim_{g\to\infty} H^*(\Kg;\bQ)^{\Mg}\cong \bQ[d_1,d_2,d_3,\cdots].
$$ 
\end{conjecture}

There has been rapid progress in the theory of the moduli space
$\overline{\mathbf M}_g$ related to the primary charactersitic 
classes $\kappa_i$ or $e_i$, namely from the celebrated works 
of Witten \cite{Wi2} and Kontsevich \cite{Ko1} to the recent developments
reaching to the Gromov--Witten invariants
(see \cite{FP, Ge} and references in them).
In contrast with this, the secondary classes $d_i$ seem to be 
beyond the reach
of any explicit study at present, except for the first one.
In relation to the first class $d_1$, it seems to be a challenging
problem to try to generalize the genus one story given in
Atiyah's paper \cite{At2} to the cases of higher genera in various ways.
We believe that the higher classes $d_i$ will also eventually
play an important role in 
hopefully deeper geometrical study of the moduli space.

\subsection{Bounded cohomology of $\Mg$}

In this subsection we consider the
mapping class group from the viewpoint of
Gromov's bounded cohomology (see \cite{Gr}).
Recall that the bounded cohomology  
(with coefficients in $\bR$) of a group $\vGa$, denoted by
$H^*_b(\vGa)$, is defined to be the cohomology of 
the subcomplex of the ordinary $\bR$--valued cochain
complex $C^*(\vGa;\bR)$ consisting of all
cochains
$$
c\co \vGa\times\cdots\times \vGa\lra \bR
$$
which are bounded as functions.
We have a natural homomorphism
$$
H^*_b(\vGa)\lra H^*(\vGa;\bR)
$$
and it reflects algebraic as well as geometric properties of
the group $\vGa$ rather closely.
For example if $\vGa$ is amenable, then $H^*_b(\vGa)$ is trivial
in positive degrees while if $\vGa$ is the fundamental group
of a closed negatively curved manifold, then the above map 
(except for degree one part)
is known to be
surjective 
by an argument due to Thurston (see \cite{Gr}).

\begin{problem}
Study the bounded cohomology of the mapping class group $\Mg$.
More precisely, determine the kernel as well as the image of the natural
map $H^*_b(\Mg)\ra H^*(\Mg;\bR)$.
\end{problem}

In particular, we may ask whether the characteristic class
$e_i\in H^{2i}(\Mg)$ can be represented by a bounded cocycle
or not.
It was remarked in \cite{Mrt3} that Gromov's general result in
\cite{Gr} implies that any odd class $e_{2i-1}$
can be represented by a bounded cocycle.
This is because, as mentioned in
\S 3, these classes are pull backs of some classical
characteristic classes of $Sp(2g,\bZ)$ which is a discrete
subgroup of $Sp(2g,\bR)$. 
It seems to be natural to conjecture
that even classes $e_{2i}$ can also be represented
by bounded cocycles.
One evidence for this was given in \cite{Mrt3} where we proved that
any surface bundle with amenable monodromy group has trivial
characteristic classes as it should be if $e_i$ were all
bounded cohomology classes. 
Meyer's signature cocycle given in \cite{Mey} is an explicit bounded cocycle
which represents $- \frac{1}{3} e_1$ (see \cite{Mrt5}) but
no other explicit bounded cocycle has been constructed for higher
odd classes.
We mention here that the cocycles for $e_i$
constructed in \cite{Mrt12, KM1} 
using trivalent graphs are far from being bounded.

\begin{problem}
Construct explicit {\it bounded} cocycles of $\Mg$ which represent
the characteristic classes $e_i$ for $i > 1$.
\end{problem}

The particular case of degree 2, namely the map
$$
H^2_b(\Mg)\lra H^2(\Mg;\bR)
$$
already deserves further investigation. Harer's determination
of $H^2(\Mg)$ in \cite{Har1} together with Meyer's result \cite{Mey} mentioned
above implies that the above map is surjective.
If $g=1$, then ${\cal M}_1=SL(2,\bZ)$ so that $H^2(SL(2,\bZ);\bR)=0$
while it is well known that $H^2_b(SL(2,\bZ))$ is infinite dimensional. 
If $g=2$, then we know that $H^2({\cal M}_2;\bR)=0$ because
${\mathbf M}_2$ is contractible by a result of Igusa \cite{Ig}.

\begin{proposition}
The two dimensional bounded cohomology $H^2_b({\cal M}_2)$ of the genus
2 mapping class group is non-trivial.
\label{prop:m2}
\end{proposition}

\begin{proof}
We consider the cochain $d_1\in C^1(\Mg;\bQ)$ described in \S 6.6.
If $g=2$, then $U=0$ so that $\delta d_1=-3\tau$ where 
$\tau\in Z^2({\cal M}_2)$ is Meyer's signature 2--cocycle for the genus 2
mapping class group. Since $H^1({\cal M}_2;\bQ)=0$, $d_1$ is
uniquely determined by the above equality.
Observe that $d_1$ is not a bounded cochain
because $d_1\co {\cal K}_2\ra \bQ$ is a non-trivial homomorphism.
We can now conclude that the bounded cohomology class 
$[\tau]\in H^2_b({\cal M}_2)$ is non-trivial.
\end{proof}

\begin{remark}
Meyer's signature cocycle has been investigated from
various points of view, see work of Y. Matsumoto \cite{Mat}
for the case of $g=2$, 
Endo \cite{En} and Morifuji \cite{Mrf2} 
for the cases of hyperelliptic mapping class groups and
Morifuji \cite{Mrf1} and Kasagawa \cite{Kas1, Kas2}
for certain geometric aspects of it.
\end{remark}

Epstein and Fujiwara \cite{EF} proved that the second
bounded cohomology of any Gromov hyperbolic group is
infinite dimensional. Although the mapping class group is not
a Gromov hyperbolic group, because it contains many free abelian
subgroup of rather high ranks, it was proved by Tromba \cite{Tr}
and Wolpert \cite{Wo1} (see also \cite{Wo2}) that the 
sectional curvature of the Teichm\"uller space with
respect to the Weil--Petersson metric is negative.

\begin{conjecture}
The mapping $H^2_b(\Mg)\ra H^2(\Mg;\bR)$ is not injective for all
$g$. More strongly, $H^2_b(\Mg)$ would be infinite dimensional.
\label{conj:bounded}
\end{conjecture}

We recall the following definition which is relevant to the above problem.

\begin{definition} (\cite{MM}).
A group $\vGa$ is said to
be {\it uniformly perfect} if there exists a natural number
$N$ such that any element $\ga\in \vGa$ can be expressed as
a product of at most $N$ commutators.
\end{definition}

Recall that $\Mg$ is known to be perfect for all $g\geq 3$ (see \cite{Har1}).

\begin{conjecture} 
The mapping class group $\Mg$ is 
{\it not} uniformly perfect for all $g\geq 3$.
\label{conj:uniform}
\end{conjecture}

It was proved in \cite{MM} that if $\vGa$ is a uniformly perfect group,
then the mapping $H^2_b(\vGa)\ra H^2(\vGa;\bR)$ is injective.
Hence if the former part of Conjecture \ref{conj:bounded} were true, then 
Conjecture \ref{conj:uniform} is also true. In the case of $g=2$, 
${\cal M}_2$ is not perfect. However it is also known that its
abelianization is finite, namely we have
$H_1({\cal M}_2)\cong \bZ/10$. 
The following result is a companion of Proposition \ref{prop:m2}.

\begin{proposition}
There is no natural number $N$ such that any element in the
commutator subgroup of the genus 2 mapping class group can
be expressed as a product of at most $N$ commutators.
\end{proposition}

\begin{proof}
If we assume the contrary, then it would follow that the value of
the cochain $d_1$ is bounded. But as was mentioned above, this
is not the case.
\end{proof}

\subsection{Representations of the mapping class group}

In this subsection, we consider various representations of the 
mapping class group.
First we mention the
Magnus representation of the Torelli group.
Let $\M1g$ be the mapping class group of $\Sg$ {\it relative}
to an embedded disk $D\subset \Sg$ as before and 
let $\bZ[\vGa]$ be the integral group ring of
$\vGa=\pi_1 (\Sg\setminus \Int D)$.
Following a general theory of so-called Magnus representation
described in Birman's book \cite{Bi},
the author defined in \cite{Mrt9} a mapping   
$$
\rho\co \M1g\lra GL(2g,\bZ[\vGa])
$$
and considered various properties of it.
It is easy to see that this mapping is injective.
However it is not a homomorphism in the usual 
sense but is rather
a crossed homomorphism. To obtain a genuine homomorphism,
we have to restrict this mapping to the
Torelli group $\I1g\subset \M1g$ and reduce the coefficients to
$\bZ[H]$ which is induced by the abelianization $\vGa\ra H$.
Then
we obtain a homomorphism
$$
\rho\co \I1g\lra GL(2g,\bZ[H])
$$
(see Corollary 5.4 of \cite{Mrt9}).

\begin{problem}
Determine whether the representaion $\rho\co \I1g\lra GL(2g,\bZ[H])$
described above is injective or not.
\end{problem}
 
Here we would like to mention that Moody proved in \cite{Mod} that the
Burau representation of the braid group $B_n$ is not faithful
for sufficiently large $n$ while it seems to be still unknown whether 
the Gassner representation of the pure braid group is faithful
or not.
\par
The augmentation ideal of $\bZ[H]$ induces a filtration of $\I1g$.
It is easy to see that this filtration is strictly coarser than
$\{\M1g(k)\}_{k\geq 1}$ which was described in \S 5. 
It would be interesting to study how they differ from each other.
\par
Besides the Magnus representation,
we have now various representations of the
mapping class group associated
to newly developed theories 
which are related to low dimensional topology. 
For example, we have projective representations of $\Mg$
arising from the conformal field thoery (see \cite{Wi1, TUY, Koh},
see also \cite{Wr, Funar} for certain explicit studies of them)
or from the theory of universal perturbative invariants of
$3$--manifolds (see \cite{LMO, Muj, MO}).
We also have Jones representations \cite{Jon} of 
the hyperelliptic mapping class group
and the Prym representations of certain subgroups of
$\Mg$ given by Looijenga \cite{Lo5}.
It seems that there are only a few results which clarify
how these representations are related to the
structure of the mapping class group.

\begin{problem}
Study various properties of the above representations of 
the mapping class group. In particular, determine the kernel as well
as the image of them.
\end{problem}

\bibliographystyle{amsplain}

\Addresses\recd
\end{document}